\documentclass[journal,twocolumn,web]{article}
\usepackage[letterpaper,top=2cm,bottom=1.9cm,left=1.9cm,right=2cm,marginparwidth=1.5cm]{geometry}
\usepackage[colorlinks=true, allcolors=blue]{hyperref}
\usepackage{amsmath,amssymb,amsfonts}
\usepackage{textcomp}
\usepackage{wrapfig}
\usepackage{overpic}
\usepackage{subcaption}

\usepackage{enumitem}
\usepackage{verbatim}
\usepackage{gensymb}
\usepackage[sorting=none, maxbibnames=99,style=ieee]{biblatex}
\usepackage[fleqn]{nccmath}
\usepackage{algorithm}
\usepackage[noend]{algpseudocode}

\algnewcommand\algorithmicinput{\textbf{INPUT: }}
\algnewcommand\Input{\item[\algorithmicinput]}
\algnewcommand\algorithmicoutput{\textbf{OUTPUT: }}
\algnewcommand\Output{\item[\algorithmicoutput]}

\MakeRobust{\Call}
\usepackage{amssymb}
\usepackage{authblk}
\makeatletter
\newcommand{\longdash}[1][2em]{%
  \makebox[#1]{$\m@th\smash-\mkern-7mu\cleaders\hbox{$\mkern-2mu\smash-\mkern-2mu$}\hfill\mkern-7mu\smash-$}}
\makeatother
\newcommand{\omitskip}{\kern-\arraycolsep}

\DeclareMathOperator*{\argmax}{argmax}

\DeclareMathOperator*{\argmin}{arg\rm{}min}

\DeclareMathOperator*{\tr}{tr}

\newcommand{\ba}{\mathbf{a}}
\newcommand{\bA}{\mathbf{A}}
\newcommand{\bB}{\mathbf{B}}
\newcommand{\bC}{\mathbf{C}}

\newcommand{\bQ}{\mathbf{Q}}
\newcommand{\bR}{\mathbf{R}}

\newcommand{\bS}{\mathbf{\Sigma}}
\newcommand{\Sel}{\mathbb{S}}
\newcommand{\bTi}{\mathbf{\Psi}} 

\newcommand{\bW}{\mathbf{W}}
\newcommand{\bx}{\mathbf{x}}
\newcommand{\by}{\mathbf{y}}

\newcommand{\bu}{\mathbf{u}}
\newcommand{\be}{\mathbf{e}}
\newcommand{\reals}{\mathbb{R}}
\newcommand{\bpi}{\mathbf{\Pi}}
\newcommand{\brr}{\mathbf{r}}

\title{Constrained optimization of sensor placement for nuclear digital twins}
\date{}
\author[1]{Niharika Karnik\footnote{nkarnik@uw.edu}}
\author[2]{Mohammad G. Abdo\footnote{Mohammad.Abdo@inl.gov}}
\author[2]{Carlos E. Estrada Perez}
\author[2]{Jun Soo Yoo}
\author[2]{Joshua J. Cogliati}
\author[2]{Richard S. Skifton}
\author[2]{Pattrick Calderoni}
\author[1]{Steven L. Brunton}
\author[1]{Krithika Manohar\footnote{kmanohar@uw.edu}}

\affil[1]{Department of Mechanical Engineering, University of Washington, Seattle, WA 98195}
\affil[2]{Idaho National Laboratory, Idaho Falls, ID 83415}

\begin{document}
\twocolumn[
  \begin{@twocolumnfalse}
\maketitle

\begin{abstract}
The deployment of extensive sensor arrays in nuclear reactors is infeasible due to challenging operating conditions and inherent spatial limitations. Strategically placing sensors within defined spatial constraints is essential for the reconstruction of reactor flow fields and the creation of nuclear digital twins. We develop a data-driven technique that incorporates constraints into an optimization framework for sensor placement, with the primary objective of minimizing reconstruction errors under noisy sensor measurements. The proposed greedy algorithm optimizes sensor locations over high-dimensional grids, adhering to user-specified constraints. We demonstrate the efficacy of optimized sensors by exhaustively computing all feasible configurations for a low-dimensional dynamical system. To validate our methodology, we apply the algorithm to the Out-of-Pile Testing and Instrumentation Transient Water Irradiation System (OPTI-TWIST) prototype capsule. This capsule is electrically heated to emulate the neutronics effect of the nuclear fuel. The TWIST prototype that will eventually be inserted in the Transient Reactor Test facility (TREAT) at the Idaho National Laboratory (INL), serves as a practical demonstration. The resulting sensor-based temperature reconstruction within OPTI-TWIST demonstrates minimized error, provides probabilistic bounds for noise-induced uncertainty, and establishes a foundation for communication between the digital twin and the experimental facility.
\end{abstract}
\end{@twocolumnfalse}
]


\section{Introduction}
\label{sec:introduction}
Safe and efficient performance of nuclear power plants requires remote monitoring, condition-based maintenance, and real-time control via data streamed from physical processes~\cite{upadhyaya2009advanced}---especially in advanced reactors (e.g., microreactors), fission batteries, small modular reactors, and integrated energy systems.
In nuclear systems, sensor capacities and real-time data streaming are severely limited, particularly for critical process responses such as coolant levels, temperature, velocity, pressure, neutron distribution and power fields. 
Furthermore, extreme operating conditions, high costs, limited accessibility and safety regulations, all impose significant constraints on sensor placement. 
The optimization of sensor placement is critical for accurate reconstruction of fields of interest, and must take into consideration not only access and safety constraints, but also the underlying physics.

\begin{figure*}[hbt!]
  \centering
    \begin{overpic}[scale=.37,percent]{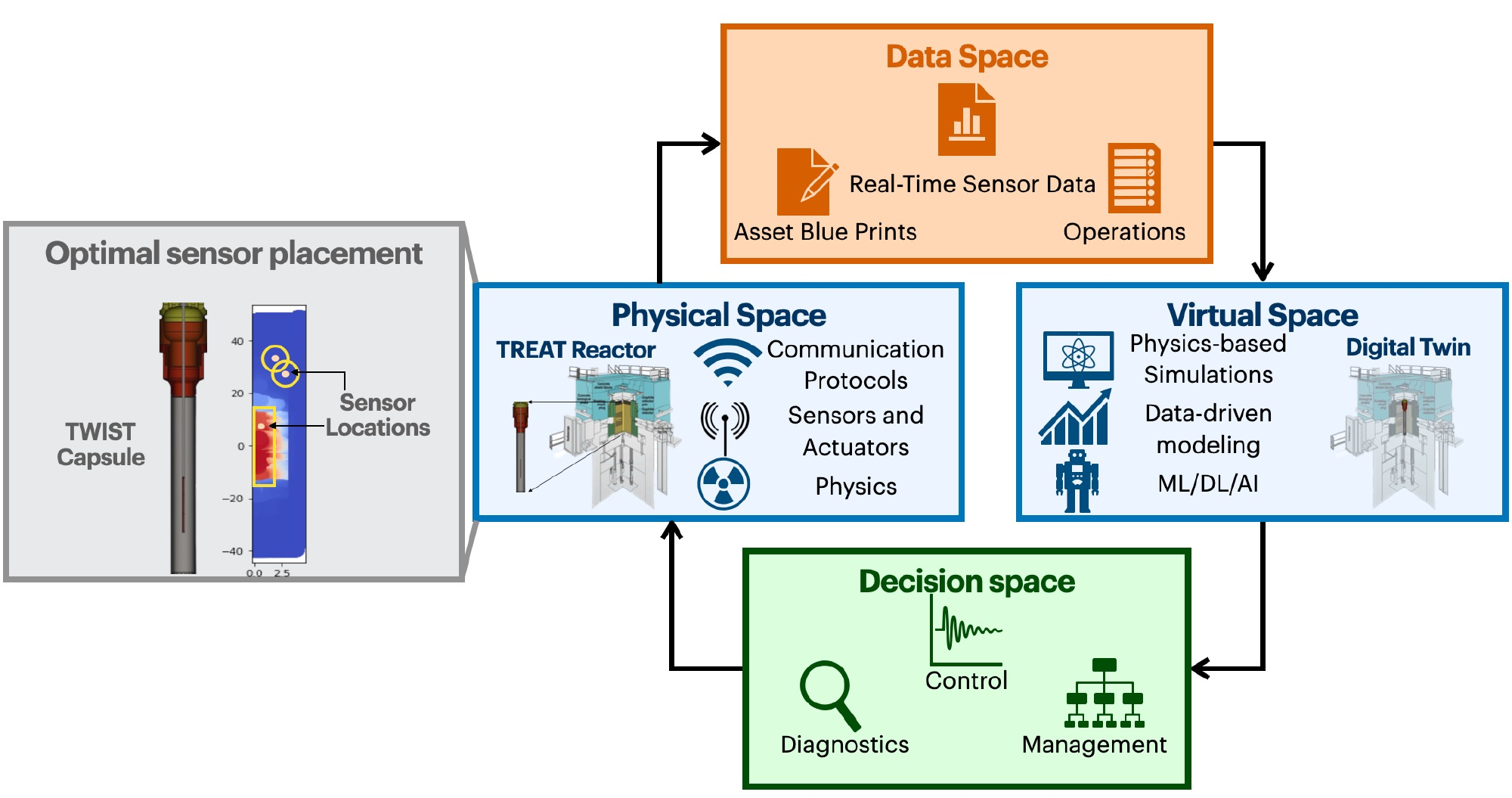}
    \put(10,33){$\by = \Sel \bx + \boldsymbol\eta$}
    \end{overpic}
  \caption{\label{fig:dt}\textbf{Digital twins in nuclear power plants.}
  Digital twin frameworks consist of a real/physical space containing physical assets, a virtual digital space containing computer-aided design (CAD) replicas, simulations and Artificial Intelligence (AI), a data space and a decision space~\cite{grieves2015digital, grieves2017digital}; all of which are enabled by sensors providing two-way communication between the virtual (ROMs and simulations) and the physical spaces. This digital twin characterizes the lifecycle of OPTI-TWIST capsule, which is inserted into the TREAT reactor at Idaho National Laboratory (INL) to test fuel compositions.}
  \vspace{-1em}
\end{figure*}

In general, sensor placement optimization is NP-hard and cannot be solved in polynomial time.  There are ${n \choose p} = n!/((n-p)!p!)$ possible combinations of choosing $p$ sensors from an $n$-dimensional state.
Exploiting low-dimensional structure inherent to the process physics is crucial for efficient optimization over complex environments. This paper explores constrained optimization of sensor placement tailored for nuclear \textit{digital twin} paradigms. Empowered by strategically placed and informative sensors, digital twins continuously stream real-time data from physical assets for digitally-enabled decision-making, control, risk assessment, and predictive maintenance, as outlined in \autoref{fig:dt}. Digital twin realizations effectively function as virtual sensors that can predict structural lifespans and ensure  structural integrity throughout product lifecycles, for example in aircrafts~\cite{tuegel2011reengineering, tuegel2012airframe}.

The goal is ``an integrated multi-physics, multi-scale, probabilistic simulation''~\cite{glaessgen2012digital, shafto2010nasa} mirroring the lifecycle of a complex system in order to reduce infeasible or costly physical testing.
In nuclear applications, digital twinning requires real-time, high-precision simulations featuring online autonomous calibration to real data via machine learning (ML) (e.g., using k-means clustering and artificial neural networks~\cite{song2022online}). Evaluating the effect of uncertainty in the digital twin, as simulated by ML models on reactor instrumentation~\cite{lin2021uncertainty}, and establishing (through sensors) a real-time two-way connection between the physical and the virtual spaces are crucial considerations for the success of digital twins~\cite{rasheed2020digital}.

Reduced-order models (ROMs) are key enablers of digital twins that compress high-fidelity, high-dimensional simulations into low-dimensional surrogate models with fewer degrees of freedom, significantly reducing computational burden while still capturing the characteristics of the relevant process~\cite{hartmann2018model,bai2005reduced}. Nuclear applications require this ability to accurately simulate high-dimensional fields with minimal computational resources and without the complexity of full-order models~\cite{wahi2011nonlinear,huang2017dimensionality}. Projection-based ROMs, which represent high-fidelity physics using low-rank/data-driven modal decompositions, have been widely adopted for modeling fluids and turbulence~\cite{sirovich1987turbulence, willcox2002balanced, noack2003hierarchy,loiseau2018sparse, rowley2005model} and nuclear core composition~\cite{abdo2019modeling, elzohery2018comparison, abdo2018analysis}. Such modal decompositions are closely related to empirical orthogonal functions for surrogate models in atmospheric sciences~\cite{franzke2005low,crommelin2004strategies}, electrodynamics~\cite{wittig2006model,koziel2012model} and heat transfer \cite{yecsilyurt1995surrogates, park1996use}; as well as model reduction of stochastic processes~\cite{everson1995karhunen} and balanced model reduction for optimal control~\cite{moore1981principal, juang1985eigensystem}. 
ROMs not only provide substantial dimensionality reduction for downstream decision-making and control, but also supply valuable physical information that can be leveraged for optimizing sensor placement.

Optimal sensor placement approaches typically optimize an objective, such as information criteria~\cite{krause2008near,joshi2008sensor}, over feasible sets of sensor configurations, framing sensor placement as a submodular selection problem~\cite{summers2015submodularity}. Such problems can be efficiently optimized for hundreds or thousands of candidate locations using convex~\cite{joshi2008sensor,brunton2016sparse} or greedy submodular optimization approaches~\cite{summers2015submodularity}.
Sensor placement in linear time-invariant systems have been optimized using gradient descent methods with similar computational complexity~\cite{chen2011h2}.
However, modern nuclear and fluid simulations have millions of grid points, making such techniques computationally intractable. 
Leveraging ROMs to optimize sensor placement, by exploiting low-dimensional patterns in data, drastically reduces the number of sensors required for accurate reconstruction of full fields~\cite{yildirim2009efficient,willcox2006unsteady,barrault2004empirical,chaturantabut2010nonlinear,drmac2016new}. Empirical ROM interpolation methods have been successfully adapted for optimizing sensor placement to minimize reconstruction error~\cite{manohar2018data,manohar2018predicting,manohar2021optimal}, including under greedy cost constraints~\cite{clark2018greedy}. However, these methods do not admit hard constraints within reactors or predict reconstruction uncertainty under a given sensor configuration with measurement noise.

This work develops a data-driven optimization approach for a nuclear component prototype, incorporating spatial constraints. In certain regions of a reactor, the placement of sensors may be constrained due to limited space availability or specific requirements dictating predetermined sensor locations, restricted areas within the reactor, fixed numbers of sensors within a region, or a minimum allowed proximity between sensors. Our target application is the Out-of-Pile Testing and Instrumentation Transient Water Irradiation System (OPTI-TWIST) prototype which is electrically heated to emulate the neutronic effect of a nuclear fuel. Production version of TWIST serves as a multi-purpose test rig for surrogate fuel rodlets, and simulate transient loss of coolant accident scenarios, to assist in qualification of an identical irradiation rig for the Idaho National Laboratory (INL) Transient Reactor Test Facility (TREAT).

We adapt data-driven methods based on modal decomposition~\cite{manohar2021optimal} to enforce these constraints during optimization, and develop placement strategies for full-field reconstruction based on sparse spatially constrained sensor measurements. Our algorithm minimizes error covariance using D-optimal design criteria, which provides an evaluation metric for a given sensor configuration and corresponding estimates of reconstruction uncertainty under noisy measurements. Using empirical and theoretical validation, the present work demonstrates the technique to be near optimal using exhaustive enumeration of all feasible sensor configurations for a low-dimensional dynamical system.
The optimized sensors under constraints are demonstrated to provide highly accurate reconstruction and uncertainty estimates under noisy measurements when compared to random placements in high-dimensional 2D heat diffusion, OPTI-TWIST steady-state and transient temperature fields with up to 40510 candidate sensor locations.
\vspace{-1.4em}
\section{Background}
This section first describes the need for optimal sensor placement in nuclear digital twinning applications through different stages within product lifecycles and product realizations. Next, we detail the reconstruction of latent flow fields from sparse sensor measurements using reduced order modeling.
\vspace{-1.4em}
\subsection{Sensing in nuclear reactors}

 A nuclear digital twin is a digital CAD replica of a physical counterpart whose complexity can vary from that of an individual fuel rodlet, heat pipe, or nuclear reactor, to that of an integrated energy system utilizing several different energy sources (e.g., wind, solar, and nuclear). Within this digital replica, visual and virtual representations of all the different sensors and actuators of nuclear components are provided~\cite{tao2018digital}. Potential application areas can leverage these digital representations for design and licensing, plant construction, training simulators, predictive operations and maintenance, autonomous operation and control, failure and degradation prediction, the generation of insights from historical plant data, and safety and reliability analyses.
 
Real-time sensor data streaming through private or Industrial Internet of Things (IoT) communication protocols is indispensable for creating digital twin architectures for nuclear applications (see \autoref{fig:dt}). The sensors provide continuous self-validation of the ML/AI models which not only reflect the current state of the dynamical system but also predict, in real time, future states of the dynamics. Current sensor technologies in the nuclear field reflect a preference that sensors be installed in easily accessible areas. Only a few algorithms have been developed for sensor placement in nuclear reactors such as the generalized empirical interpolation method \cite{argaud2018sensor}, reinforcement learning \cite{gu2022radiation}, and a directed graph approach for minimizing postulated faults with maximum imperceptibility \cite{upadhyaya2011optimal}. 

Optimal sensor technologies empower digital twins by critically enabling the integration of in-field and real-time raw data into the virtual replica of a physical prototype at any point during its product 
lifecycle \cite{darvishi2020sensor}, which include the design, manufacturing, service and retirement stages. Our data-driven sensor placement methodology is designed for optimizing sensors in the design stage. Virtual sensors are converted into physical sensors and validated based on experimentation throughout design, manufacturing and service. If new constraints arise during production, optimal sensing techniques can incorporate them in real time and then suggest the next best set of optimal sensor locations. 
\vspace{-1.3em}
\subsection{Sparse sensing for reconstruction}


The core of our work is the reconstruction of latent fields $\bx\in\reals^n$ from $p$ noise-corrupted sensor measurements $\by\in\reals^p$ 
\begin{equation}
    \by =  \Sel \bx + \boldsymbol\eta, 
\end{equation}
where $\boldsymbol{\eta}$ consists of zero-mean, Gaussian independent and identically distributed (i.i.d.) components, and $\Sel \in \reals^{p\times n}$ is the desired sensor (measurement) selection operator. In nuclear applications, the number of measurements $p$ is severely limited relative to the large dimensionality of the latent field. We encode the field dynamics as a linear combination of spatial basis modes $\boldsymbol{\psi}_k(\xi)$ weighted by time-varying coefficients
\begin{equation}\label{eqn:pod}
    \bx(t) \approx \sum_{k=1}^{r} a_k(t) \boldsymbol{\psi}_k(\xi).
\end{equation}
For each field or full state $\bx$ at a fixed $t$, the vector $\ba$ composed of the $r$ coefficients $a_k(t)$ defines a low-rank embedding of the form
$$\bx = \bTi_r\ba,$$
where the modes $\boldsymbol{\psi}_k$ comprise the  columns of $\bTi_r$. This basis, which can be built from spectral or data-driven decomposition methods, is typically chosen so that the embedding dimension is as small as possible, i.e., $r\ll n$.

Given this assumption, high-dimensional states can be directly recovered from measurements via the maximum likelihood estimate of the basis coefficients, $\hat\ba = (\Sel \bTi_r)^\dagger\by$:
\begin{equation}\label{eqn:gappy}
    \hat\bx = \bTi_r(\Sel \bTi_r)^\dagger\by,
\end{equation}
known as gappy proper orthogonal decomposition (POD)~\cite{everson1995karhunen}. The gappy estimator is well-posed when the number of sensors equals or exceeds $r$. Importantly, the inherent compressibility of physical fields enables a drastic reduction in the number of sensors required for high-fidelity reconstruction. The critical enabler for \emph{sparse sensing} is the fact that nuclear processes are strictly governed by a small set of underlying physics. As we shall see, strategic selection of sensor measurements---based on noisy flow physics---allows for an extremely small number of deployed sensors to be used.
\vspace{-1em}
\subsection{Data-driven basis}
The data embedding rank dictates the minimum number of sensors required for reconstruction, necessitating a choice of basis with the lowest possible rank. Given full state data sampled from physics/CFD simulations $\mathbf{X} = \begin{bmatrix}
    \bx_1 &\dots &\bx_m
\end{bmatrix}$, the proper orthogonal decomposition (POD)~\cite{berkooz1993proper} provides the minimal rank approximation to data
\begin{equation}\label{eqn:podopt}
    \argmin_{\bTi,\operatorname{rank}(\bTi)=r} \| \mathbf{X} - \bTi\bTi^T\mathbf{X}\|_F^2,
\end{equation}
where the low-rank embedding is given by the projection of the data onto orthogonal POD modes, $\bTi_r^T\mathbf{X}$. The solution to ~\eqref{eqn:podopt} is computed using the singular value decomposition of the data matrix, $\mathbf{X} = \mathbf{UD} \mathbf{V^*}$, where the leading $r$ left singular vectors comprise the desired POD modes 
\begin{equation}\label{eqn:svd}
    \bTi_r = \mathbf{U}_r = \begin{bmatrix}
        \bu_1 &\bu_2 &\dots &\bu_r
    \end{bmatrix}.
\end{equation}
The singular values (diagonal entries of $\mathbf D$), quantify the decreasing energy contribution of each successive mode and determine the truncation rank. Most physical data have much fewer degrees of freedom than the ambient data dimension, allowing a very small choice of $r$. The cumulative energy captured by the leading $r$ modes is $\sum_{i=1}^r d_i/\sum_i d_i$. In practice, the smallest possible $r$ capturing 90-99\% of  cumulative energy above the noise threshold is used as the model truncation rank. Therefore, POD is also the workhorse of projection-based model order reduction, used for projecting governing equation terms onto POD modes to obtain highly computationally expedient surrogate models for high-fidelity physics.

\section{Methodology}
This section describes our methodology for optimizing sensors for reconstruction in constrained settings. First, we describe the ROMs used to set up a linear reconstruction problem for recovering high-dimensional fields from sparse measurements (i.e., sensors). We then characterize our sensor placement optimization objective in terms of reconstruction error. Next, we detail how our greedy algorithm selects the next optimal sensor in terms of a strategic projection operator, and modify the selection step to enforce the necessary constraints while still maintaining optimality.

\subsection{Sensor placement for reconstruction}
The placement of sensors is defined by a measurement selection operator $\Sel \in \reals^{p \times n}$ that optimally recovers modal mixture $\ba$ from sensor measurements $\by$. This measurement selection operator $\Sel$ encodes point measurements with unit entries in a sparse matrix
\begin{equation}
    \Sel = \begin{bmatrix} \be_{\gamma_1} & \be_{\gamma_2} & \hdots & \be_{\gamma_ p}\end{bmatrix}^T,
\end{equation}
where $\be_j$ are canonical basis vectors for $\mathbb{R}^n$, with a unit entry in component $j$ (where a sensor should be placed) and zeroes elsewhere. Here, $\gamma = \{\gamma_1, \gamma_2, \hdots, \gamma_p\} \subset \{1, 2, \hdots, n\}$ denotes the index set of sensor locations with cardinality $p$. 
Sensor selection then corresponds to the components of $\mathbf{x}$ that were chosen to be measured:
\begin{equation}
   \mathbf{\Sel x} = \begin{bmatrix} x_{\gamma_1} & x_{\gamma_2} & \hdots & x_{\gamma_p}\end{bmatrix}^T.
\end{equation}
The selection of sensors is based on the optimal estimation of the entire state vector $\mathbf{x} \in \mathbb{R}^n$ from $p$ experiment outputs $\mathbf{y} \in \mathbb{R}^p$ with additive i.i.d. Gaussian noise $\eta_i \sim \mathcal{N} (\mathbf{0},\beta^2)$ in each measurement $\by_i$:
\begin{equation}
\label{eqn:measurement}
    \by = \Sel\bTi_r\ba + \boldsymbol{\eta}.
\end{equation}
The values of $\mathbf{x}$ at unmeasured locations can be recovered by solving a linear system of equations for the basis coefficients via the Moore-Penrose pseudoinverse of $\Sel\bTi_r$ (gappy POD~\eqref{eqn:gappy}):
\begin{equation*}
    \hat\bx = \bTi_r(\Sel\bTi_r)^\dagger\by.
\end{equation*}
The row indices of $\bTi_r$ correspond to sensor locations in the state space that effectively condition the matrix inversion, enabling accurate reconstruction of the estimated state $\mathbf{\hat{x}}$. 

Optimal design of experiments~\cite{fisher1936design} for estimation problems involves the strategic selection of a set of experiments to gather sufficient information about the domain, enabling accurate predictions for measurements where experiments were not performed.  Statistical criteria, such as A, D and E-optimality, are used to select the set which minimizes or maximizes different properties of the Fisher information matrix. Fisher information~\cite{fisher1936design} measures the amount of information a random variable contains about the estimated parameter, such as its true mean or standard deviation. The Fisher information matrix defines covariance matrices associated with maximum-likelihood estimates and is $(\Sel \bTi_r)^T(\Sel \bTi_r)$ in our case. A-optimal designs minimize the trace of the inverse of the Fisher information matrix, whereas E-optimal designs maximize the minimum eigenvalue of the information matrix. D-optimal designs ~\cite{smith1918standard} minimize the generalized variance of the parameter estimates by maximizing the determinant of the Fisher information matrix ~\cite{sengupta2004generalized}. 

Optimal design for gappy estimation involves placing sensors at limited points in the domain to accurately reconstruct flow fields over the entire domain. In contrast to classical optimal design in which each sensor can be used multiple times out of a set of candidate sensors, candidate sensors can only be used once in the gappy framework. In this setting, design of experiments aims to optimize the sensor selection $\Sel$ to optimize statistics of the estimation error $\ba-\hat\ba$, an $r$-dimensional random variable with zero mean and covariance
\begin{equation}
\label{eqn:Co}
    \boldsymbol{\Sigma} = \text{Var}(\mathbf{a} - \mathbf{\hat{a}}) = \beta^2((\Sel \bTi_r)^T(\Sel \bTi_r))^{-1}.
\end{equation}
 The eigenvalues of this covariance matrix characterize the statistical and geometric measures of estimation error ``size"~\cite{friendly2013elliptical}, shown in \autoref{demo-table}.
\begin{table}[hbt!]
\begin{center}
\begin{tabular}{|c| c|c|} 

 \hline
 \footnotesize{Measure} & \footnotesize{Formula} & \footnotesize{Geometry}\\ 
 \hline\hline
 \footnotesize{Generalized variance} & \footnotesize{$\det(\boldsymbol{\Sigma}) = \Pi_i \lambda_i$} & \footnotesize{area, (hyper)volume} \\ 
 \hline
 \footnotesize{Average variance}  & \footnotesize{$\text{tr}(\boldsymbol{\Sigma}) = \sum_{i} \lambda_i$} & \footnotesize{linear sum}\\
 \hline
 \footnotesize{Maximal variance}  &  \footnotesize{$\lambda_{\max}$} & \footnotesize{maximum dispersion}\\

 \hline

\end{tabular}
\caption{\label{demo-table} Statistical and geometric measures for error covariance~\cite{friendly2013elliptical}}
\end{center}
\end{table}

Generalized variance, defined by $\det(\boldsymbol{\Sigma})$, characterizes correlations among pairs of variables. When it is large, the variables have little correlation with each other; when it is small, the variables are strongly correlated. On the other hand, average variance, given by $\tr(\boldsymbol{\Sigma})$, is  the sum of the population variances. A-optimal criteria minimize this average variance, while E-optimal criteria minimize the maximal variance of $\boldsymbol{\Sigma}$.
The variance, which measures the uncertainty in the estimated response, should be small for minimal deviation between estimated and true values~\cite{smith1918standard}. 

We consider D-optimal design for flow field reconstruction  with information matrix $(\Sel \bTi_r)^T(\Sel \bTi_r)$, which depends on the selected sensors $\Sel(\gamma)$. The determinant objective maximizes the information volume via maximization of its determinant, given a budget of $p$ sensors. The maximizing sensor set of this criterion is also the maximizer of its logarithm
\begin{equation}
\label{eqn:detobj}
    \gamma_* = \argmax_{\gamma, |\gamma|= p} \; \log \det((\Sel \bTi_r)^T(\Sel \bTi_r)). 
\end{equation}

When $p = r$,~\autoref{eqn:detobj} is equivalent to the maximizer of $\log |\det(\Sel \bTi_r)|$.
Direct optimization of this criterion leads to a brute force combinatorial search. This sensor placement approach builds upon the empirical interpolation method (EIM) \cite{barrault2004empirical} and discrete empirical interpolation method (DEIM) \cite{chaturantabut2010nonlinear} which select the best interpolation points for evaluating nonlinear terms in projection-based reduced order models. However, these methods do not directly optimize statistics of the error or minimize error covariance. In the next section, we develop a greedy strategy for optimizing sensor selection under constraints built upon the pivoted QR factorization~\cite{drmac2016new,manohar2018data, manohar2018predicting, manohar2021optimal}, and analyze the reconstruction performance with respect to D-optimal design criteria. 

\subsection{Column-pivoted QR decomposition with spatial constraints}
\label{ssec:Constrained QR}
The QR factorization with column pivoting decomposes a matrix $\bW\in \reals^{m \times n}$ into a unitary matrix $\bQ$, an upper-triangular matrix $\bR$, and a column permutation matrix $\bpi$, such that $\bW \bpi = \bQ \bR$.
As described above, each column index of $\bTi_r^T$ corresponds to a single sensor location in the state space. We applied QR pivoting to the transpose of our basis, i.e. $\bW =\bTi_r^T$, and use the permutation matrix to store information about the sensors selected. 
The pivoted QR decomposition is a greedy algorithm for optimizing \autoref{eqn:detobj} that, in each iteration) selects a new column pivot (sensor location) with maximal two-norm, then subtracting from every other column vector its orthogonal projection onto the pivot column \cite{manohar2018data, drmac2016new}. This projection is given by a Householder reflector that maps any vector $\boldsymbol{\nu}$ to $-\operatorname{sign}\nu_1\sigma\mathbf{e}_1$, where $\sigma = \|\boldsymbol{\nu}\|_2$ and $\nu_1$ is the first component of $\boldsymbol{\nu}$
\begin{equation}
\label{eqn:householder}
    \mathbf{H}(\boldsymbol{\nu}) = \mathbf{I} -  \frac{(\boldsymbol{\nu} + \text{sign}(\nu_1)\sigma \mathbf{e}_1)(\boldsymbol{\nu} + \text{sign}(\nu_1)\sigma \mathbf{e}_1)^{T}}{\sigma(\sigma + |\nu_1|)}.
\end{equation}
 Householder projections effectively zero out the subdiagonal components of column vectors in each iteration to induce upper-triangular structure in $\bW$, constructing $\bR$ in place. Householder reflectors can be written in the standard form $\mathbf{I} - 2\mathbf{u}\mathbf{u}^T$, where $\bu$ has unit norm.

\begin{figure*}[t!]
  \centering
  \includegraphics[width=0.6\textwidth]{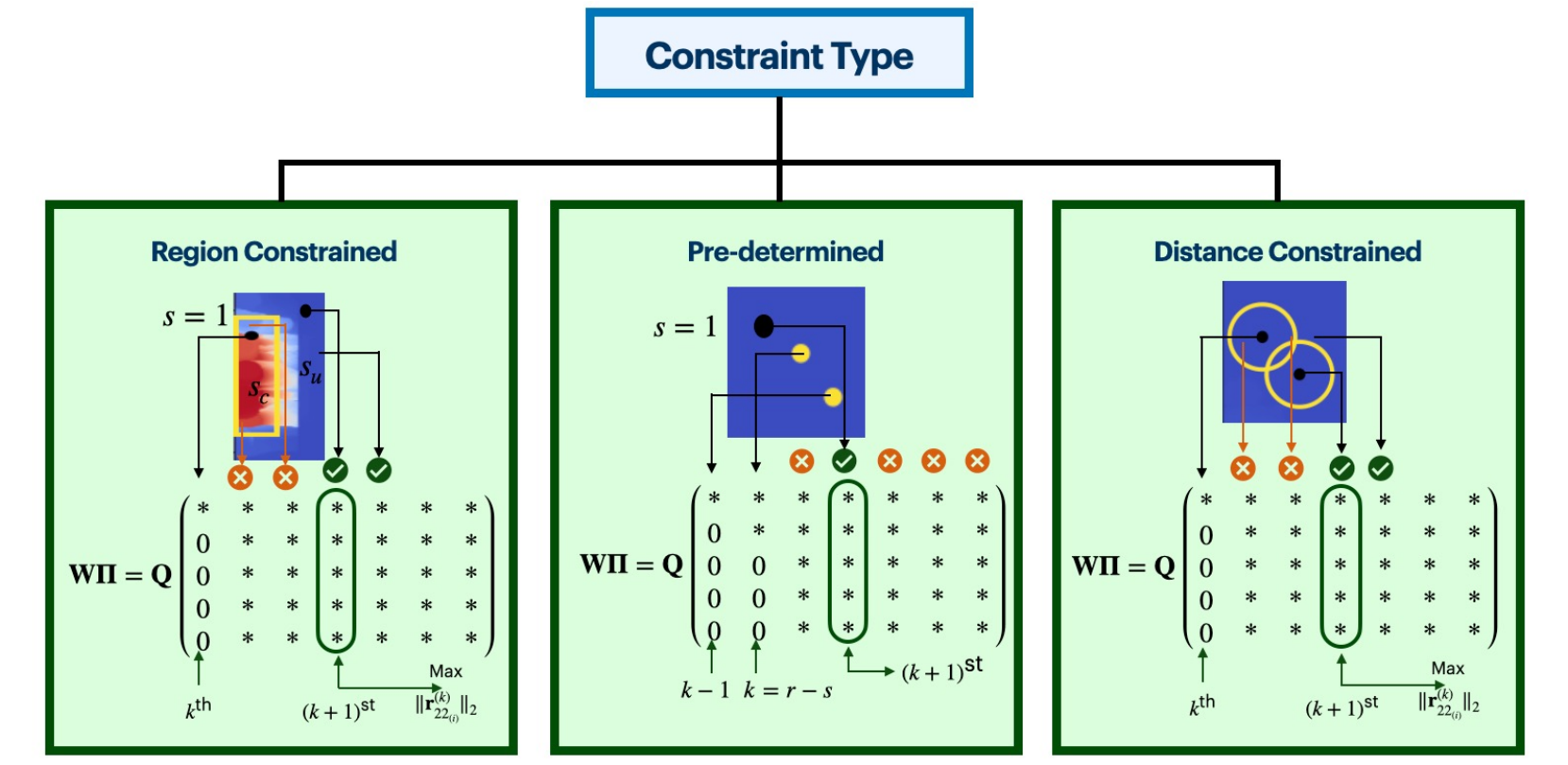}
  \caption{\label{fig:flowchart} {\bf Greedy selection of the next sensor} involves choosing the next pivot column of $\bTi^T$ from the set of allowable sensor locations specified by the constraint.}
\end{figure*}

Consider the following partial QR factorization at the $k^{\text{th}}$ iteration in the pivoting procedure: 
\begin{equation}
    \bW\bpi = \bQ \bR = \bQ 
    \begin{bmatrix} \bR_{11}^{(k)} & \bR_{12}^{(k)} \\ {0} & \bR_{22}^{(k)} 
    \end{bmatrix},
\end{equation}
where $\bQ \in \reals^{m\times m}$ is orthogonal, $\bR_{11}^{(k)} \in \reals^{k \times k}$ is upper triangular, $\bR_{12}^{(k)} \in \reals^{k \times (n-k)}$, $\bR_{22}^{(k)} \in \reals^{(m-k) \times (n-k)}$, and $\bpi \in \reals^{n \times n}$ is the permutation matrix containing information about the first $k$ chosen sensors \cite{businger1965linear, clark2020multi, gu1996efficient}.
In unconstrained QR pivoting, the $(k+1)^{\text{st}}$ iteration selects a column from the submatrix $\bR_{22}^{(k)}$ with the maximal two-norm, then swaps the selected column with the $(k+1)^{\text{st}}$ column while updating permutation indices
\begin{equation}
    l = \argmax_{i = 1, \hdots , n-k} \|\brr_{22_{(i)}}^{(k)}\|_2 .
\end{equation}

Constraints are integrated within this step, by forcing the pivot column index to be selected from the latest set of allowable indices based on the constraints under consideration. The $k + 1^{\text{st}}$ iteration in constrained optimization selects the pivot column with largest 2 norm, $\brr_{22_{(i)}}^{(k)}$, from the constrained/unconstrained set of allowable column indices in $\bR_{22}^{(k)}$.
The QR pivoting algorithm results in the following diagonal dominance structure in $\bR$:
\begin{equation}
\label{eqn:diagDom}
    |\bR_{ii}|^2 \geq \sum_{j = 1}^{k}|\bR_{jk}|^2, \, \, 1 \leq i \leq k \leq r.
\end{equation}
Constraints are imposed in the final $r-s$ steps of pivoting, ensuring that the largest contributing terms in the objective function expansion remained unaffected:
\begin{fleqn}
    \begin{equation}
        \begin{aligned}[b]
            &\log \det ((\Sel\boldsymbol{\Psi_r})^T\Sel\boldsymbol{\Psi_r})    
                =\log \left(\prod_{i=1}^r \bR_{ii}^2\right)\\\nonumber
            & = \log \bR_{11}^2 + \log \bR_{22}^2 + \dots + \log\bR_{(r-s)(r-s)}^2 + \dots + \log \bR_{rr}^2.
        \end{aligned}
    \end{equation}
\end{fleqn}
The main driver of this optimization is the point at which constraints are introduced into the pivoting procedure, as allowing upper triangularization to proceed normally in the starting iterations maximizes the leading diagonal entries of $\bR$, ensuring that domain-specific constraints do not drastically affect the diagonal dominance property, but only the trailing $\bR_{ii}$, which are optimized by choosing the best pivot from the allowable locations. 
The three types of spatial constraints handled by the algorithm (\autoref{fig:flowchart}) are:

\begin{enumerate}[topsep=0pt,itemsep=0em]
    \item \textbf{Region constrained: } This type of constraint arises when we can place either a \textit{maximum} of or \textit{exactly} $s$ sensors in a certain region, while the remaining $r-s$ sensors must be placed outside the constraint region.
    \begin{itemize}
        \item{\textbf{Maximum:}} This case deals with applications in which the number of sensors in the constraint region should be less than or equal to $s$. In each iteration a pivot column (sensor location) is chosen from the set of all columns until $s$ selected pivots lie in the constrained region. Successive pivots with the largest 2-norm are selected from among the unconstrained column indices. 
        \item{\textbf{Exact:}} This case deals with applications in which the number of sensors in the constraint region should equal $s$. The algorithm follows the same procedure as the maximum sensor placement case if the number of sensors in the constraint region equals or exceeds $s$. However, if there are fewer than $s$ sensors in the constrained region, the algorithm forces the deficit of sensors to be placed in the constraint region at the end of the pivoting procedure. 
     \end{itemize}
    \item \textbf{Predetermined:} This type of constraint occurs when a certain number of sensor locations are already specified, and  optimized locations for the remaining sensors are desired. The strategy employed selects pivots from among all column indices of $\bR_{22}^{(k)}$ until the iterate $k $ equals $r-s$, then imposes the selection of  user-specified sensor locations in the final $s$ iterations.
    \item \textbf{Distance constrained:} This constraint enforces a minimum distance $d$ between selected sensors. Accordingly, the first pivot is the column index of $\bTi_r^T$ with maximal 2-norm, the default (unconstrained) base step.  The $(k+1)^{\text{st}}$ iterate now selects the pivot column with maximal 2-norm from among the remaining columns of $\bR_{22}^{(k)}$ that are at least distance $d$ away from the previous $k$ selected sensors. This is an \emph{adaptive} constraint because the set of allowable sensor indices is updated with each pivoting iteration to remove the $d$-neighborhood of the $k$th sensor.
\end{enumerate}
\newcommand\sForAll[2]{ \ForAll{#1}#2\EndFor} 
\newcommand\sIf[2]{ \If{#1}#2\EndIf}          
\begin{algorithm}[b!]
    \caption{QR Pivoting with Adaptive Constraints }\label{algorithm1}
    \begin{algorithmic}[1]
  \State \textbf{Input:} $\bW$ = Input matrix, $s_c$ = Constrained/Predetermined sensor indices, $s$ = indices allowed in the constrained region/ Total predetermined sensors, {\footnotesize CONSTRAINTS} = Type of Constraint, $d$ = Distance from previous sensors
  \State \textbf{Output:} Sensor indices
  \State \textbf{Function Call:} $\gamma$ = {\footnotesize CONSTRAINED}QR($\Psi_r^T$, {\footnotesize CONSTRAINTS}, $s_c$=[], $s$=[], $d$=2)
    \Procedure{constrainedQR($\bW$, {\footnotesize Constraints},$s_c$,$s$,$d$)}{}
        \State $r,n \gets$\Call{Size}{$\bW$}  \Comment{$r$ = number of desired sensors}
        \State $\bR \gets$\Call{Copy}{$\bW$}
        \State $\bQ \gets eye(m)$
        \State $\gamma \gets [1,2,\dots,n]$   \Comment{sensor index set}
        \sForAll{$k\in 1,....,r$}{
            \State $dlens\gets$\Call{{\footnotesize {\footnotesize computeConstraints}}}{{$dlens$,$\gamma$,{\footnotesize CONSTRAINTS},
            $k$,$s_c$,$s$,$d$}}
            \State $l \gets \argmax dlens$
            \State $ \boldsymbol{\nu} \gets \bR[k:m,k-1+l]$ 
            \State \Call{swap}{$\bR[k:m,k],\bR[k:m,k-1+l]$}
            \State \Call{swap}{$\gamma_k,\gamma_{k-1+l}$} \Comment{update $k$th sensor}
            \State $u \gets \Call{Householder}{\boldsymbol{\nu}}$ \Comment{compute $\mathbf{H}(\boldsymbol{\nu})$} 
            \State $\bR[k:m,k:n] \gets (I -  2uu^T)\bR[k:m,k:n]$  
            \State $\bQ[:,k:n]\gets (I - 2uu^T) \bQ[:,k:n]$
        }
    \State \Return $\gamma_{1:r}$
\EndProcedure
\end{algorithmic}   
\end{algorithm}

Although we mainly consider the minimal allowable number of sensors to be $p=r$, the truncation rank of the basis, additional sensors can be added for redundancy and robustness through the oversampling optimization proposed by B. Peherstorfer et al~\cite{peherstorfer2020stability}. 
\begin{algorithm}[t!]
    \centering
    \caption{Subroutine for Constraints}\label{algorithm2}
    \begin{algorithmic}[1]
  \State \textbf{Input:} $s_c, s, r, d$, {\footnotesize CONSTRAINTS} (see \autoref{algorithm1}), $dlens$ = $\bR$ column norms, $\gamma$ = sensor index set, $k$ = current sensor index
  \State \textbf{Output:} $dlens$\Procedure{{\footnotesize computeConstraint}($dlens$,$\gamma$,{\footnotesize Constraints},$k$,$s_c$,
  $s$,$r$,$d$)}{}
            \If{{\footnotesize CONSTRAINTS} = RegionConstrainedMax}
                \If{$|\gamma| \in s_c < s$} 
                    \State$dlens$ = $dlens$
                \Else 
                \State $dlens[s_c] = 0$
                \EndIf
    
            \ElsIf{{\footnotesize CONSTRAINTS} = RegionConstrainedExact}
                \If{$|\gamma| \in s_c < s$} 
                    \State$dlens$ = $dlens$
                        \If{$r > k \geq (r - (s-|\gamma|)$} 
                            \State$dlens[! s_c] = 0$
                        \EndIf
                \Else \Comment{Number of sensors in constrained region $ > s$ }
                    \State $dlens[s_c] = 0$
                \EndIf
            \ElsIf{{\footnotesize CONSTRAINTS} = Pre-determined}
                \If{$k = (r - s)$} 
                    \State$dlens[!s_c]$ = 0 
                \EndIf
            \Else{ {\footnotesize CONSTRAINTS} = DistanceConstrained}
                \If{$p \in dlens $ \textbf{and} $\|p -\gamma[:]\|_{2} \leq d$} 
                    \State $dlens[p] = 0$
                \EndIf
            \EndIf
\EndProcedure
\end{algorithmic}
\end{algorithm}

\subsection{Uncertainty analysis}
\label{sec:UQ}
Under noisy measurements, errors in estimation are transmitted into reconstruction errors. Geometrically, estimation errors are characterized by (hyper-)ellipsoids in $r$ dimensions whose axes describe these errors. Statistically, the confidence intervals for the estimation of states are characterized by the \emph{$\eta$-confidence ellipsoid} that contains $\mathbf{a} - \mathbf{\hat{a}}$ with probability $\eta$
\begin{equation}
    E_{\alpha} = \{ \mathbf{z} | \mathbf{z}^T \boldsymbol{\Sigma}^{-1} \mathbf{z} \leq \alpha\},
\end{equation}
where $\mathbf{\Sigma}$ is the covariance matrix in \autoref{eqn:Co} and $\alpha = F_{\chi_n^2}^{-1} (\eta)$ is the cumulative distribution function of a $\chi$-squared random variable with $r$ degrees of freedom. An important scalar measure of the quality of estimation is the volume of this ellipsoid%
\begin{equation}
        \begin{aligned}[b]
            &\textbf{vol}(\epsilon_{\alpha}) = \frac{(\alpha \pi)^{r/2}}{\Gamma (\frac{r}{(2+1)})} \det \boldsymbol{\Sigma}^{1/2} \\
            & =\frac{(\alpha \pi)^{r/2}}{\Gamma (\frac{r}{(2+1)})} \det ((\beta^2((\Sel\bTi_r)^T(\Sel\bTi_r))^{-1})^{1/2}),
        \end{aligned}
    \end{equation}
where $\Gamma$ is the gamma function. D-optimal designs minimize the volume of the ellipsoid, which is inversely proportional to the determinant of our information matrix. A small volume for the $\eta$-confidence ellipsoid implies a strong correlation between the estimation errors in each component. Under Gaussian noise assumptions, $3\sigma$ standard deviations computed from the diagonal entries of the covariance matrix $\Sigma_{ii}$ measure the uncertainty in predicting the $i^{th}$ component, establishing error bounds for reconstructing flows from noisy measurements. 

The following text, for ease of notation, uses $\Theta = \Sel\bTi_r$ to signify the mode measurement matrix, hence $\bS=\Theta^T \Theta$ represents the Fisher Information. In order to quantify the uncertainty in each reconstructed component of the state under noisy measurements, we analyze the expected covariance of the full state fluctuations. Klishin et al~\cite{klishin2023data} provide the following estimate of the expected state covariance for a regularized gappy estimator, which now depends on the covariance of our measurement model as follows 
\begin{equation}
    \begin{aligned}
    \label{eqn:state_uncertainty}
    \mathbb{E}[\Delta\hat{\bx} \Delta\hat{\bx}^T] &=\Psi_r \frac{\bS^{-1}}{\beta^2}\Theta^T\mathbb{E}[\Delta y\Delta y^T]\Theta \frac{\bS^{-1}}{\beta^2} \Psi_r^T,
    \end{aligned}
\end{equation}
where $\mathbb{E}[\Delta y\Delta y^T] = \beta^2\mathbf{I}$ for uncorrelated noise. The standard deviation in the reconstruction of each grid-point can be calculated through the diagonal entries of this matrix, providing an uncertainty heatmap of the whole reconstructed domain based on the sensor configuration $\Sel$~\cite{klishin2023data}
\begin{equation}
\begin{aligned}
\mathbf{F} & \equiv \Psi_r \bS^{-1} \frac{\Theta^T}{\beta^2} \\
\sigma_i & =\beta \sqrt{\sum_j\left(\mathbf{F}_{i j}\right)^2}.
\end{aligned}
\end{equation}
Furthermore, model recalibration for digital twinning can be informed by analyzing the distribution of each component of the estimated coefficients $\hat{\ba}$, even when the true coefficients are unavailable. The predicted mean $\mu_i=\mathbb{E}[\hat{\ba}_i]$ is estimated by averaging measurements over the training data 
\begin{equation}
    \label{eqn:mean_ahat}
    \mathbb{E}[\hat{\ba}_i] = \bS^{-1} \frac{\Theta^T}{\beta^2} \mathbb{E}[\by_i]
\end{equation}
where $\by_i$ are each component of noisy measurements with standard deviation $\beta$.
Similarly, we can estimate the expected covariance in components of $\hat{\ba}$ using the diagonal entries of
\begin{equation}
\begin{aligned}[b]
    \label{eqn:avg_covariance}
    \mathbb{E}[\Delta\hat{\ba} \Delta\hat{\ba}^T] = \frac{\bS^{-1}}{\beta^2}\Theta^T\mathbb{E}[\Delta y\Delta y^T]\Theta \frac{\bS^{-1}}{\beta^2}.
\end{aligned}
\end{equation}
The diagonal part $\mathbf{T}$ of the covariance matrix can be used to calculate the standard deviation in the distribution of the estimated POD coefficients. 
\begin{equation}
\begin{aligned}
  \label{eqn:Diagonal_cov}
    \mathbf{T} & \equiv \bS^{-1} \frac{\Theta^T}{\beta^2} \\
\sigma_i & =\beta \sqrt{\sum_j\left(\mathbf{T}_{i j}\right)^2}.
\end{aligned}
\end{equation}
The predicted standard deviation and mean of the POD coefficients together provide statistical metrics to measure uncertainty in the reconstruction of the flow field due to noise measurements when true readings are unavailable. 
In nuclear digital twins these error bounds can be used to detect the divergence of sensor readings from expected values. Statistics of the error provide means to ``flag'' or ``signal'' re-calibration of the digital twin, detect anomalies and classify erroneous readings.
\vspace{-1em}
\section{Results}
This section demonstrates the constrained and unconstrained sensor placement algorithm on a randomly generated state space system, the 2D heat diffusion through a thin plate, and the OPTI-TWIST prototype. In the randomized system, all possible sensor placements given a fixed budget of sensors are computed to demonstrate the near optimality of our approach. 
Next, we investigate reconstruction of temperature fields in 2D heat diffusion with a constant heat source, a simplified model of the OPTI-TWIST heater. Uncertainty analysis is conducted on noisy measurements for varying signal-to-noise ratio (SNR). Using constrained optimized sensor placement with our approach, we reconstruct, with minimal error, the flow field inside OPTI-TWIST---in comparison to randomly selected sensor locations.  
\vspace{-1em}
\subsection{Discrete random state space}
We first demonstrate the near optimality of constrained QR pivoting by using a low-dimensional linear time-invariant (LTI) system. The dimensionality of this system is small enough to enumerate all possible placements in order to empirically compare the constrained QR placements with the optimum placements. Consider the following LTI system:
\begin{subequations}
\begin{align}
    \dot{\bx} &= \bA\bx + \bB\bu & \mathbf{x} \in \mathbb{R}^n, \mathbf{u} \in \mathbb{R}^q\\           
    \by &= \bC\bx & \mathbf{y} \in \mathbb{R}^p
\end{align}
\end{subequations}
with randomly generated system $\bA$, measurement $\bC$, 
\begin{figure}[hbt!]
	\centering
	\begin{subfigure}{.44\textwidth}\label{unconBF}
		\includegraphics[width=\textwidth]{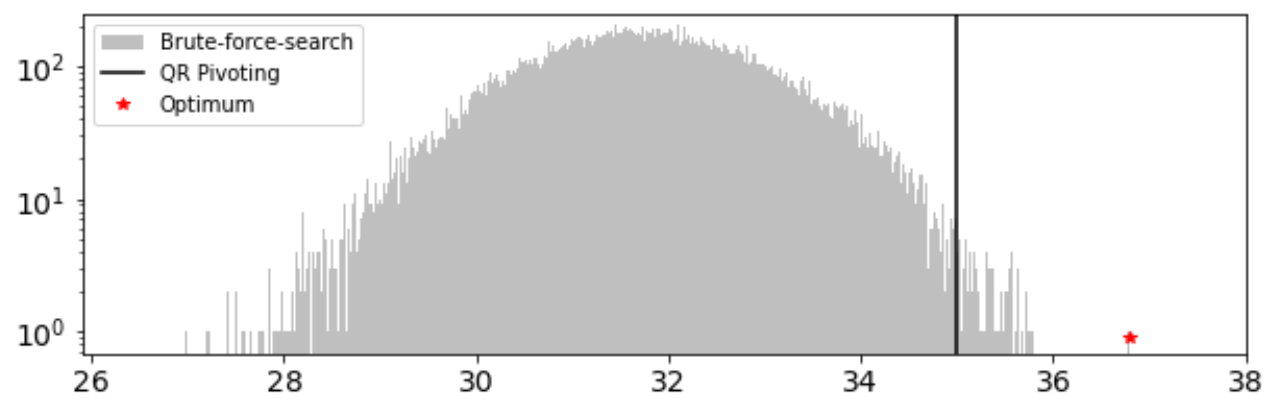}
		\caption{\label{fig:rssunconst}Unconstrained QR pivoting}
	\end{subfigure}
	\begin{subfigure}{.44\textwidth}
		\includegraphics[width=\textwidth]{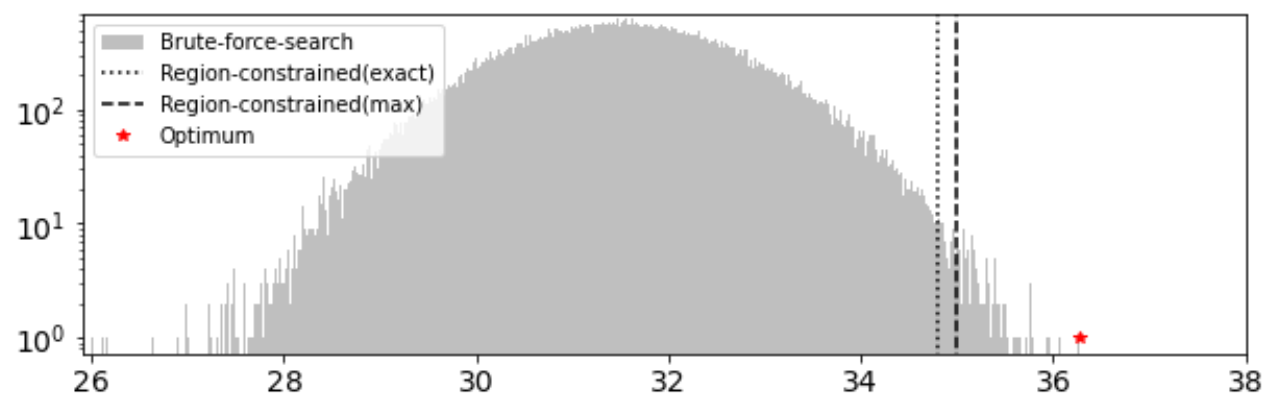}
		\caption{\label{fig:regconst}Region-constrained QR pivoting}
	\end{subfigure}
	\begin{subfigure}{.45\textwidth}
     \begin{overpic}[scale=.38,percent]{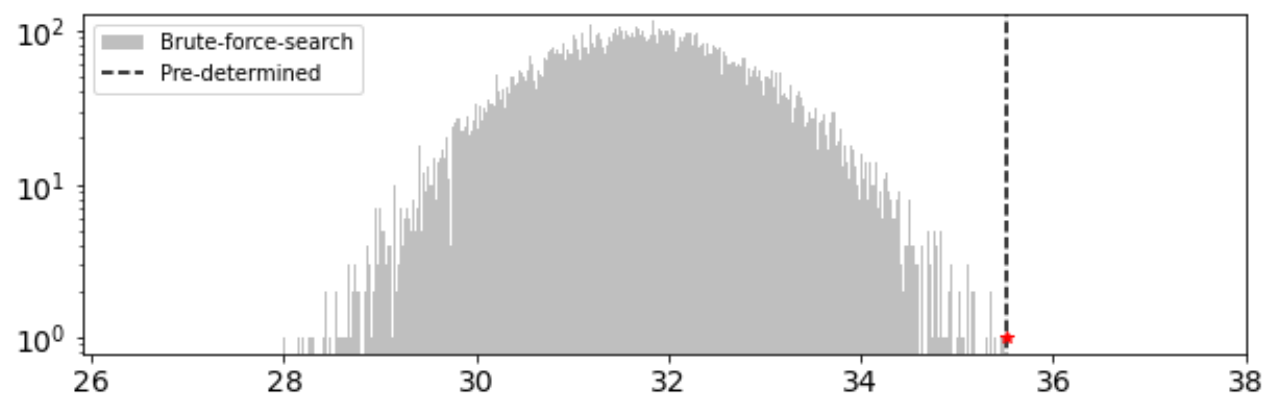}
        \put(-3,38){\rotatebox{90} {Number of Placements}}
    \end{overpic}
		\caption{\label{fig:preconst}User-specified QR pivoting}
	\end{subfigure}
	\caption{\textbf{Enumeration of $\log \det(\Sel\boldsymbol{\Psi_r})^T(\Sel\boldsymbol{\Psi_r})$ (X-axis) over all possible placements of 7 out of 25 candidate locations (100,000-500,000 possible placements binned into histograms).} The introduction of constraints into QR optimization results in a log determinant that is near optimal (optimum shown in red) for the three types of constraints.}
\end{figure}%
and actuation $\bB$ matrices sampled i.i.d. from a normal distribution, $n = 25$ states and $p = q = 25$ randomized measurements and actuators, respectively. $\bx$, $\bu$, and $\mathbf{y}$ represent the state, input, and output vectors, respectively. We empirically studied the optimality of our proposed algorithm by computing all possible $r = 7$ sensor selections, leading to a brute force search of ${n \choose r} = 480,700 $ choices of $\Sel$. The log determinant of $\Sel\boldsymbol{\Psi_r}$ was evaluated for all possible combinations of the seven sensors, then binned into the histogram shown in \autoref{fig:rssunconst}. This computation is only tractable for a small state dimension---as even for $n=100$, the brute-force search results in $O(10^{10})$ complexity. The optimization outcome (determined via \autoref{eqn:detobj}) for sensors selected using the QR pivoting approach is represented by the solid black line in \autoref{fig:rssunconst}. This sensor selection is observed to be nearly D-optimal, exceeding 99.74\% of all candidate placements.

We studied region-constrained pivoting by allowing only $s\le 2$ sensors to be selected from the first $s_c=5$ components of the state (the constraint region). A brute-force search across all possible combinations of $s$ sensors in the constraint region (and $r-s$ elsewhere) was carried out, resulting in $\frac{(n-s_c)!}{(n-s_c-(r-s))!(r-s)!} \frac{s_c!}{(s_c-s)!s!} = 155,040$ possible combinations in selecting the seven sensors binned in \autoref{fig:regconst}. \autoref{fig:regconst} compared our first strategy (i.e., placing exactly $s = 2$ sensors in the constraint region in the first two iterations of pivoting (dotted line)) with another strategy in which a maximum of two sensors were allowed in the constraint region throughout the pivoting procedure (dashed line). Our exact approach has a log determinant exceeding 99.78\% of all possible region-constrained sensor placement options, while the max approach is observed to exceed 99.87\% of all possible combinations. Thus, both approaches provide a near-optimal subset of region-constrained sensors for reconstruction. 

In predetermined sensor placement, a specified number $s$ of sensors were selected in advance, leaving the algorithm to optimize those that remain. Our strategy runs unconstrained QR pivoting to select the first $r-s$ pivots (sensors), then selects the predetermined sensors in the remaining $s$ pivoting iterations. The results of the log determinant objective evaluated for our optimization strategy (dashed line) are compared against a brute force search across all possible candidate placements that contain the two predetermined sensors reflected in \autoref{fig:preconst}. Our strategy outperforms 99.99\% of all possible placements, exhibiting near-optimal solutions.
\begin{figure}[hbt!]
  \centering
  \includegraphics[width=0.42\textwidth]{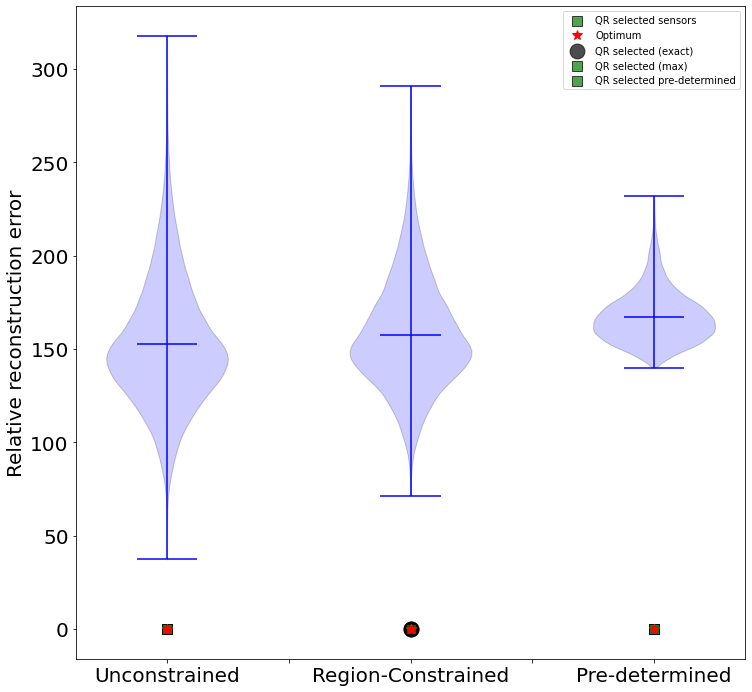}
  \caption{\label{fig:RSS-Ensemble} \textbf{Reconstruction error comparison.} Proximity between the brute-force optimum and QR selected sensors for unconstrained, region-constrained, and predetermined sensor placement leads to orders of magnitude lower reconstruction error ($\epsilon \sim \mathcal{O}(10^{-15})$) compared to random placements. Incorporating constraints results in accuracy comparable to that of the optimal placement (red stars).}
\end{figure}
These results show that the introduction of constraints results in minimal distance to the true optimum. We analyze the negligible effect of this distance on the reconstruction error
\begin{equation}\label{eqn:relrecerror}
    \text{\footnotesize{Relative reconstruction error}}\footnotesize{(\epsilon) = \frac{\|\bx -  \bTi_r(\Sel \bTi_r)^\dagger\by\|_{2}}{\|\bx\|_2} \times 100.}
\end{equation}
We compared the reconstruction achieved via each set of constrained QR sensors with sensor placements sampled from the mean of the $\log\det$ distributions. (These represent the most likely sub-optimal sets to be chosen at random.) The reconstruction error for each of these randomly placed sensors was then calculated (see the blue violin plots in \autoref{fig:RSS-Ensemble}, where the green square and circle reflect the reconstruction error of the QR-selected sensor locations for the different constraint cases). Random sensor placements with a sub-optimal $\log\det$ objective fall between 31 and 32, resulting in a highly inaccurate relative reconstruction error $(\epsilon)$ between 30 and 350. The QR-optimized strategy for unconstrained/constrained sensor placement results in significantly lower reconstruction error $\epsilon \sim \mathcal{O}(10^{-15})$:
\begin{figure*}[hbt!]
	\centering
 \begin{subfigure}{.2\textwidth}
		\includegraphics[width=\textwidth]{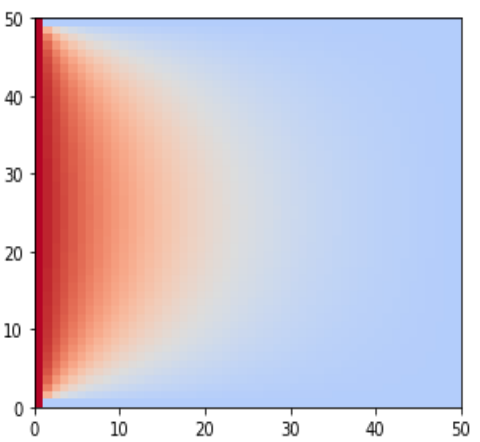}
		\caption{\label{fig:heatplate3}\footnotesize{Ground Truth}}
	\end{subfigure}
	\begin{subfigure}{.19\textwidth}
		\includegraphics[width=\textwidth]{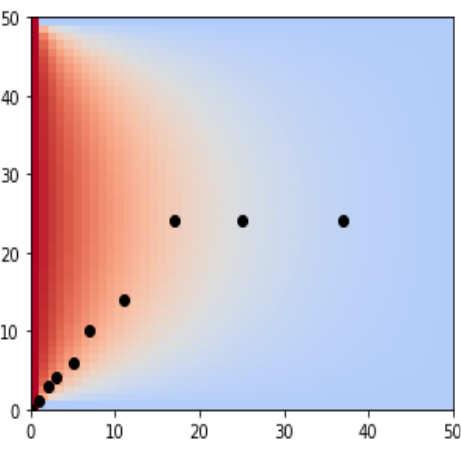}
		\caption{\label{fig:heatplate3}\footnotesize{QR $\epsilon$ = 0.0064}}
	\end{subfigure}
 \begin{subfigure}{.19\textwidth}
		\includegraphics[width=\textwidth]{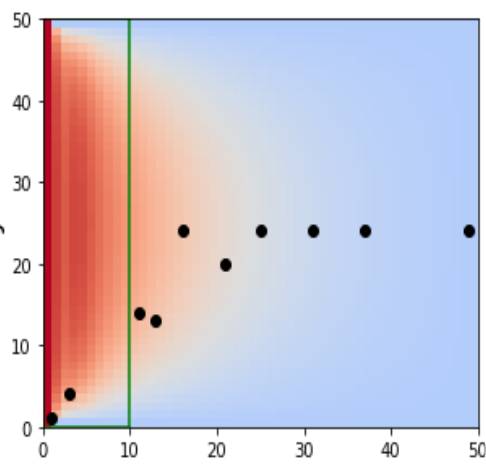}
		\caption{\label{fig:heatplate4}\footnotesize{Exact,Max $\epsilon$ = 1.5085}}
	\end{subfigure}
 \begin{subfigure}{.25\textwidth}
		\includegraphics[width=\textwidth]{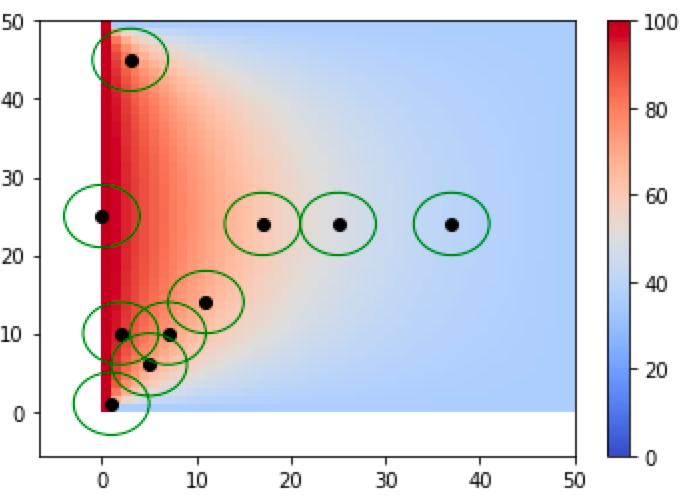}
		\caption{\label{fig:heatplate8} \footnotesize{Distance $\epsilon$ = 0.0103}}
	\end{subfigure}
 \begin{subfigure}{.2\textwidth}
		\includegraphics[width=\textwidth]{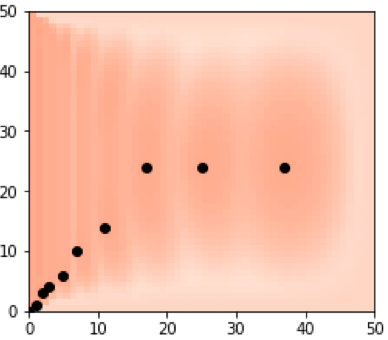}
		\caption{\label{fig:heatunc1}\footnotesize{QR}}
	\end{subfigure}
 \begin{subfigure}{.205\textwidth}
		\includegraphics[width=\textwidth]{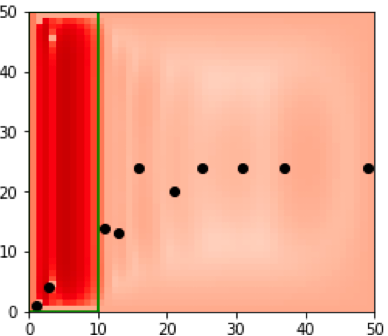}
		\caption{\label{fig:heatunc2}\footnotesize{Exact and Max equal}}
	\end{subfigure}
 \begin{subfigure}{.24\textwidth}
		\includegraphics[width=\textwidth]{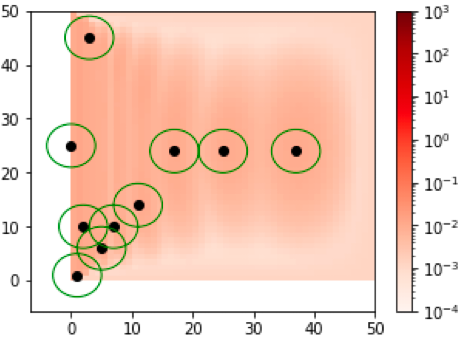}
		\caption{\label{fig:heatplate7}\footnotesize{Distance}}
	\end{subfigure}
	\caption{\label{fig:heat}\textbf{Reconstruction of the temperature field through selected sensors along with uncertainty in reconstruction.} Uncertainty heatmaps (e,f,g) correspond to placements/reconstructions (b,c,d) respectively. Reconstruction of the temperature field at $t = 1000$, based on the different constraints demonstrate that constraining sensors far away from the heater region result in higher reconstruction error (c,d) and higher uncertainty (f,g) than unconstrained optimization (b,e) respectively, which favors sensors adjacent to the heat source.}
\end{figure*}
\begin{figure}[hb!]
\begin{subfigure}{.2\textwidth}
		\includegraphics[width=\textwidth]{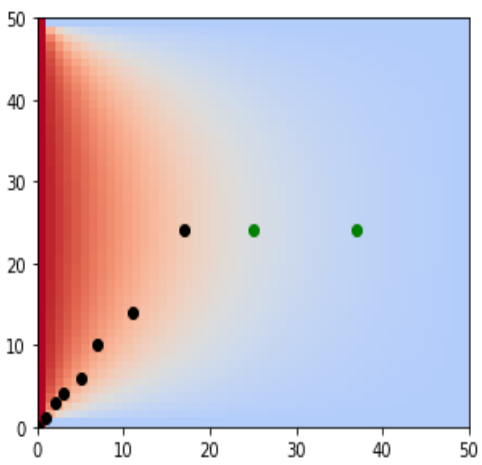}
		\caption{\label{fig:heatplate6}Pre $\epsilon$ = 0.0127}
	\end{subfigure}
	\begin{subfigure}{.26\textwidth}
		\includegraphics[width=\textwidth]{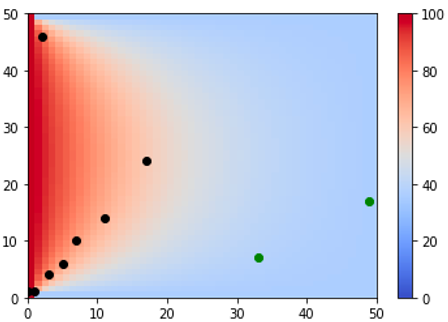}
		\caption{\label{fig:heatplate7}Pre $\epsilon$ = 0.0111}
	\end{subfigure}
\begin{subfigure}{.22\textwidth}
		\includegraphics[width=\textwidth]{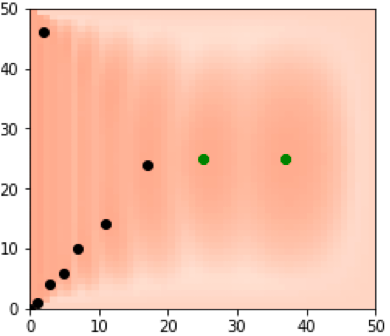}
		\caption{\label{fig:heatunc4}Pre}
	\end{subfigure}
\begin{subfigure}{.28\textwidth}
		\includegraphics[width=\textwidth]{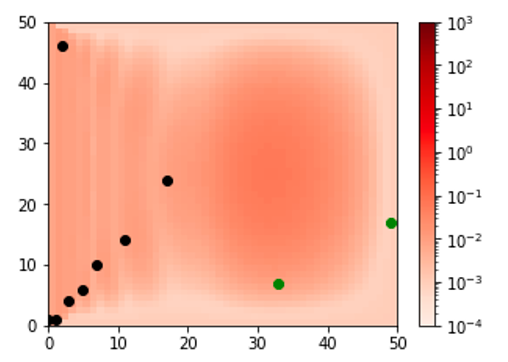}
		\caption{\label{fig:heatunc5}Pre }
	\end{subfigure}
	\caption{\label{fig:heat_pre}\textbf{Informed sensor placement lowers uncertainty in reconstruction. } Two different pre-determined sensor layouts (a,b-shown in green) may lead to similar reconstruction errors (a,b) but increase reconstruction uncertainty~(d) when sensors are distant from QR-optimal locations (c).}
\end{figure}
\begin{figure*}[hbt!]
\centering
\begin{subfigure}{0.49\textwidth}
    \includegraphics[width=\textwidth]{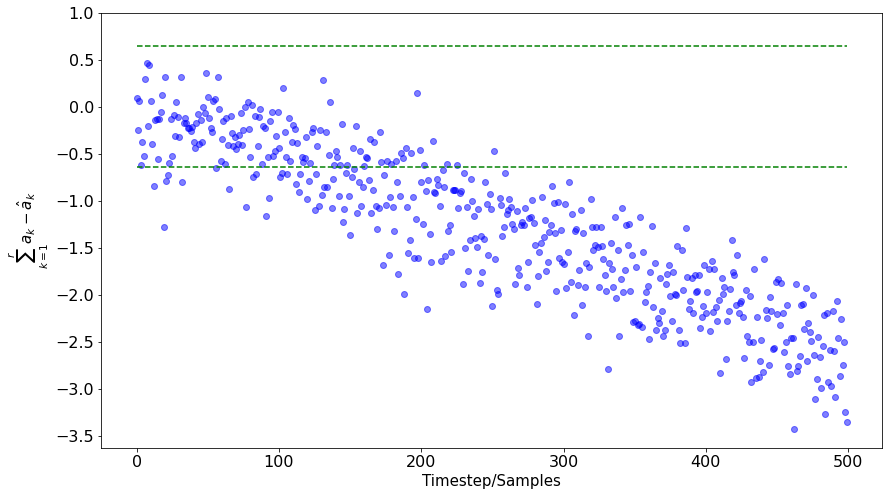}
	       \caption{\label{fig:heat_r10model} \footnotesize{Rank 10 model}}
\end{subfigure}
\begin{subfigure}{0.49\textwidth}
\includegraphics[width=\textwidth]{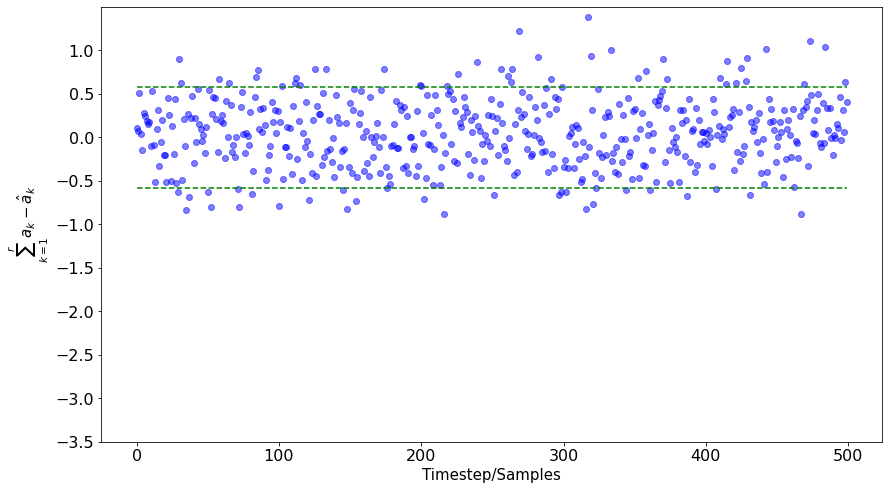}
	       \caption{\label{fig:heat_r20model}\footnotesize{Rank 20 model}}
	\end{subfigure}
 
 \begin{subfigure}{0.48\textwidth}
		\includegraphics[width=\textwidth]{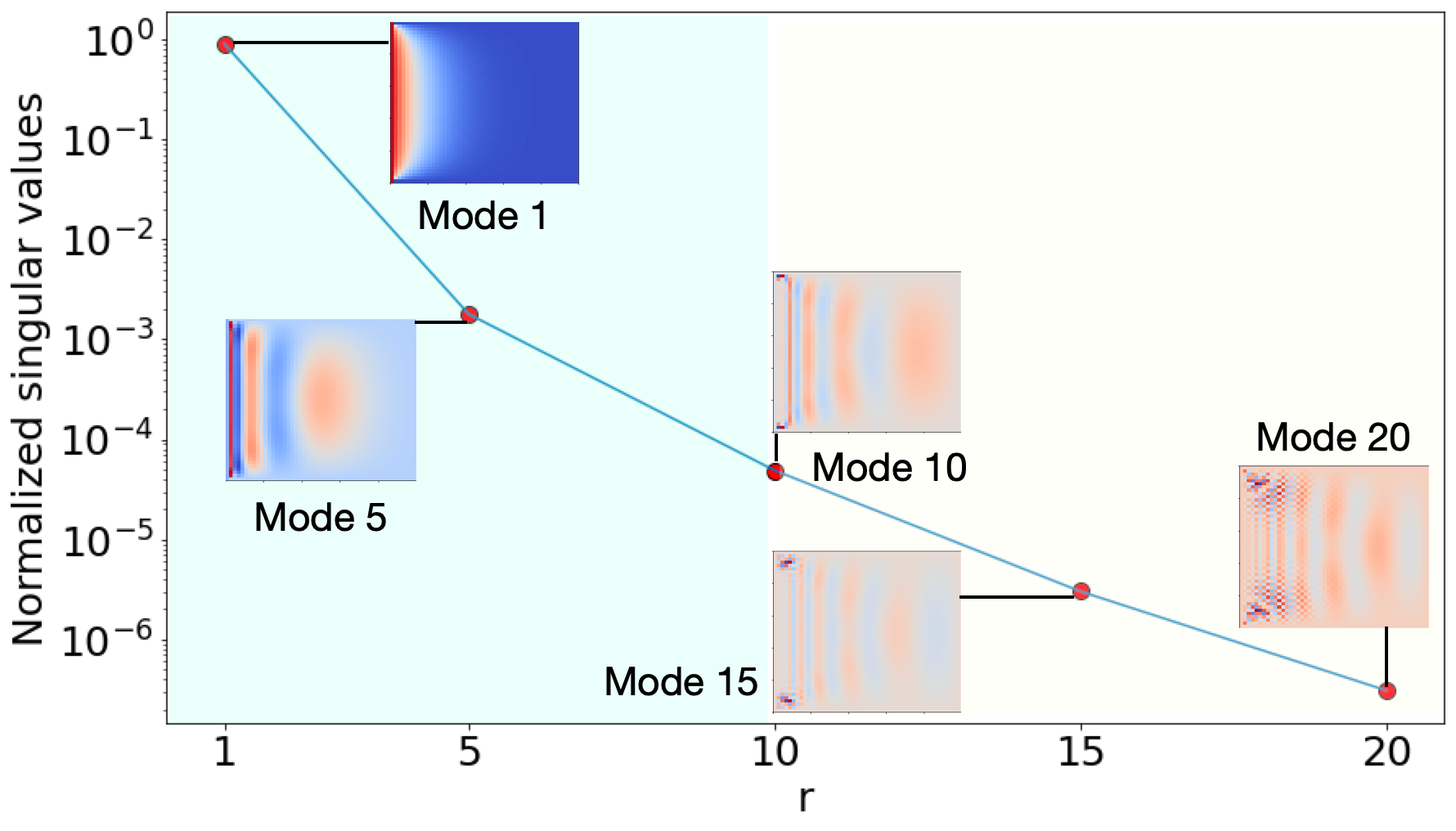}
	       \caption{\label{fig:heatPOD}\footnotesize{Leading POD modes by energy contribution.}}
        \end{subfigure}
 \begin{subfigure}{0.5\textwidth}
 \includegraphics[width=\textwidth]{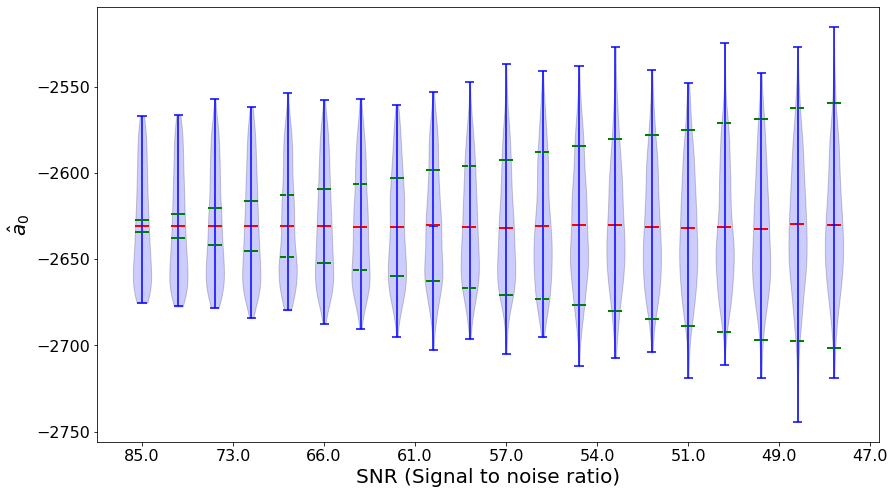}
		\caption{\hspace{-0.6em}\label{fig:boundssnr} \footnotesize{Increasing uncertainty in estimation of $\hat{\ba}_k$ with a decreasing SNR.}}
	\end{subfigure}
\caption{\label{fig:heat_modes_and_noise}\textbf{Analysis of uncertainty in estimation due to noisy sensor data.} Uncertainty estimation reveals rank 10 models are not sufficiently descriptive of dynamics after t =100s (a) under noisy measurements. A higher rank 20 approximation is required (b) despite the leading 10 POD modes capturing 99\% of the energy (c). Predicted statistics of the estimated coefficients, such as the standard deviation $3\sigma$ (green) and mean (red) effectively bound the estimated $\hat{\ba}_0$ under increasing noise.}
\end{figure*}

We conclude that the proximity between the brute-force optimum and our greedy placements produces negligible loss in reconstruction performance. Vastly sub-optimal random placements illustrate the inverse relationship between the objective function and performance: the lower the log determinant, the higher the reconstruction error. However, this system was randomly generated, and the dynamics do not evolve according to localized features in state space that allow one sensor to capture information on  spatially correlated states. Next, we demonstrate the algorithm on a heat diffusion model that allows for spatial and physical interpretation of the resulting sensors.
\vspace{-1.3em}
\subsection{2D heat flow through a thin plate}
Temperature fields are reconstructed via constrained/unconstrained sensor placement in a 2D plate undergoing thermal diffusion from a heat source, based on a simplified model of the OPTI-TWIST diffusion. Temperatures at the boundaries are fixed at $100 \degree C$ at $x=0$, and $36 \degree C$ elsewhere in $\partial D$. The initial temperature throughout the domain $D$ at $t=0$ is also $36 \degree C$. Heat transfer from the heat source over the domain is governed by the heat equation  
\begin{align*}
    &u_t = \alpha(u_{xx} + u_{yy})  &x,y\in D\\
    &u(x,y,t) = 36 &x,y\in \partial D\setminus x=0 \\
    &u(0,y,t) = 100 &
\end{align*}
where $u$ is the desired temperature and $\alpha=2 \, (\text{mm}^2/\text{s}$) is the thermal diffusivity constant. The solution of the partial differential equation is simulated for 1000 time steps using finite differences with $\Delta x,\Delta y = 1$ and $\Delta t= 0.125$. We reconstruct the temperature fields and analyze the uncertainty in reconstruction of each pixel due to adding i.i.d. Gaussian noise $\boldsymbol{\eta} \sim \mathcal{N} (0,0.01)$ to the measurements for the unconstrained, constrained (\autoref{fig:heat}), and predetermined optimized sensor placements (\autoref{fig:heat_pre}). A total of $r = 10$ sensors are selected for reconstruction, and $s = 2$ sensors are allowed in the constrained region or are predetermined. Distance constraints impose a Euclidean distance of at least 2 between selected sensors.

Similar to the nuclear OPTI-TWIST (\autoref{ssec:OPTI-TWIST}), optimized placements favor sensors near the heat source. Constraining sensors  away from the heat source results in higher reconstruction errors and uncertainty than unconstrained optimization. In this example, both region-constrained max and exact ($s=2$ within $x<10$) optimization result in identical sensor placements (\autoref{fig:heatplate4}), with only two sensors near the heat source. This results in higher error because of high-energy modal contributions adjacent to the heater~(\autoref{fig:heatPOD}). Removing heater adjacent sensors results in higher uncertainty of approximately 10\degree C near the heater~ (\autoref{fig:heatunc2}). The distance-constrained sensor placements, which also placed six sensors near the heat source, 
result in the lowest reconstruction error under constraints (\autoref{fig:heatplate8}). When predetermined sensor locations are close to the unconstrained optimal locations as in \autoref{fig:heatplate6}, the reconstruction errors and uncertainty are low (\autoref{fig:heatunc4}). Optimized sensor locations (unconstrained) are placed along the propagating wavefront, as in the two rightmost sensors. Fixing predetermined sensor locations away from this wavefront results in a noticeable increase in reconstruction uncertainty in the right half of the domain~(\autoref{fig:heatunc5}), and strengthens the case for data-driven placement strategies. The optimized sensors reflect the highest energy amplification in the POD modes, which occur near the heat source (\autoref{fig:heatPOD}). The leading POD modes capture this diffusion of heat from the heat source boundary to the rest of the domain (\autoref{fig:heatPOD}). Approximately 96\% of the cumulative energy is captured by the leading two POD modes, while the remainder capture only 4\%. 

Gappy POD was used to estimate the rank 10 and rank 20 model coefficients from noisy measurements (test dataset of 500 snapshots), and compare estimation errors with our predicted uncertainty analysis to test the descriptive capability of the different rank truncations. The standard deviations $\sigma_i$ computed from the diagonal entries of the covariance matrix $\bS_{ii}$ measures the uncertainty in predicting the $i$th component. Approximately 498.5 out of 500 reconstruction errors should lie within $3\sigma_i$ of the mean. As more modes are included in the POD approximation, the selected sensors capture more information about the underlying physics of heat diffusion. As seen in \autoref{fig:heat_modes_and_noise}, the rank 10 POD model fails to capture the underlying physics after time interval $t = 100$, whereas the rank 20 model is more descriptive of the dynamics over a longer time interval $t = 500$. With more modes and sensors, the error covariance in each component narrows and the $3\sigma$ bounds become tighter.
\,
When clean measurements or ground truth coefficients are unavailable, bounds on distributions of estimated $\hat{\ba}$ are useful to inform recalibration of digital twins. The uncertainty estimation for any POD coefficient, for example $\hat{\ba}_0$, can be bounded  using the predicted standard deviation or 3$\sigma$ (\autoref{eqn:Diagonal_cov}, shown in green in \autoref{fig:boundssnr}) and mean (\autoref{eqn:mean_ahat}, shown in red in \autoref{fig:boundssnr}). As the SNR decreases, the dynamics of heat diffusion are corrupted by noise, resulting in wider distributions of the estimated state. Note that predicted uncertainty bounds are more accurate at higher noise levels due to numerical rank approximation error overwhelming the low noise contribution to error. In other words, uncertainty analysis is more accurate under larger ratios of measurement noise to model error.

\begin{figure*}[hbt!]
 \hspace{0.5em}
    \begin{subfigure}{.13\textwidth}
        \begin{overpic}[scale=.67,percent]{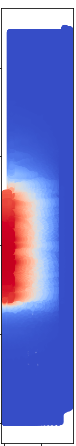}
          \put(-6.5,34){\rotatebox{90}{Height [cm]}}
    \end{overpic}
		
		\caption{\label{fig:TWISTgroundTruth} \scriptsize{Ground Truth}}
	\end{subfigure}
	\begin{subfigure}{.1\linewidth}
            \includegraphics[width=0.7\textwidth]{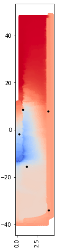}
		\caption{\label{fig:TWISTrand}\scriptsize{Random}}
	 \end{subfigure}
  \hspace{1.5em}
	\begin{subfigure}{.075\textwidth}
		\includegraphics[width=\textwidth]{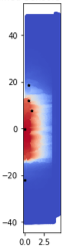}
		\caption{\label{fig:TWIST1}\scriptsize{QR}}
	\end{subfigure}
 \hspace{1.5em}
	\begin{subfigure}{.155\textwidth}
		\includegraphics[width=\textwidth]{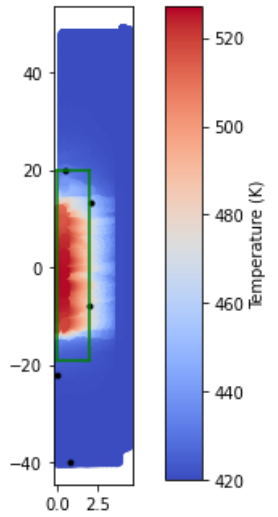}
		\caption{\label{fig:TWIST3}\scriptsize{Constrained} }
	\end{subfigure}
  \hspace{2em}
 \begin{subfigure}{.1\linewidth}
		\includegraphics[width=0.7\textwidth]{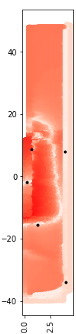}
		\caption{\label{fig:TWIST_unc_rand} \scriptsize{Random}}
	\end{subfigure}
  \hspace{1.5em}
\begin{subfigure}{.075\textwidth}
		\includegraphics[width=\textwidth]{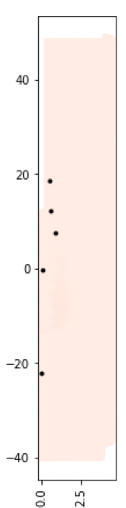}
		\caption{\label{fig:TWIST_unc_unconst} \scriptsize{QR}}
	\end{subfigure}
  \hspace{1.5em}
\begin{subfigure}{.125\textwidth}
		\includegraphics[width=\textwidth]{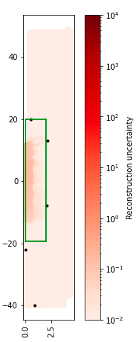}
		\caption{\label{fig:TWIST_unc_const}\footnotesize{Constrained} }
	\end{subfigure}
\caption{\label{fig:TWISTprofile}\textbf{OPTI-TWIST temperature profile reconstructions through different sensor layouts and uncertainty in estimation caused by noisy sensor measurements.} Unconstrained optimization places sensors near the heater region (c), resulting in highly accurate reconstruction with $\epsilon = 0.168$ (a), with constrained optimized sensors resulting in comparably high accuracy $\epsilon = 0.174$ (d). Random sensor placement (b) results in inaccurate reconstructions ($\epsilon = 25.24 $) and large estimation uncertainty (e) compared to that of optimized sensor locations (f,g).}
\end{figure*}

In summary, when uncertainty in sensor measurement is known, this framework enables estimation of the expected error distribution as a function of measurement noise, as well as study of the growth in error as the sensor noise increases. When sensor uncertainty is unknown, filtering and Bayesian inference techniques may be used for uncertainty quantification with these linear embeddings. 
The developed algorithm can handle reconstruction of flow fields in the presence of constraints and noise with high accuracy. Uncertainty analysis of predicted states plays a key role in detecting erroneous readings in digital twins. In the next example we reconstruct temperature flow field for a gravitational advection driven physical system, OPTI-TWIST.

\vspace{-1em}
\subsection{Steady-state simulation of the OPTI-TWIST prototype} 
\label{ssec:OPTI-TWIST}

The next example follows the new design paradigm suggested by digital twins. In traditional design practice, modeling and simulation insights are often leveraged at the experimental design stage in order to build physical models and place sensors. However, limitations regarding space, installation, cost, and signal fidelity of the experimental device pose challenges in deploying the desired number of sensors. Our holistic approach integrates experimental constraints, Computational Fluid Dynamics (CFD) simulations, and optimization objectives (reconstruction) in a principled way to optimize the placement of sensors in the design phase of the digital twin.

\begin{figure}[hbt!]
    \begin{subfigure}{.24\textwidth}
        \centering
		\includegraphics[width=\textwidth]{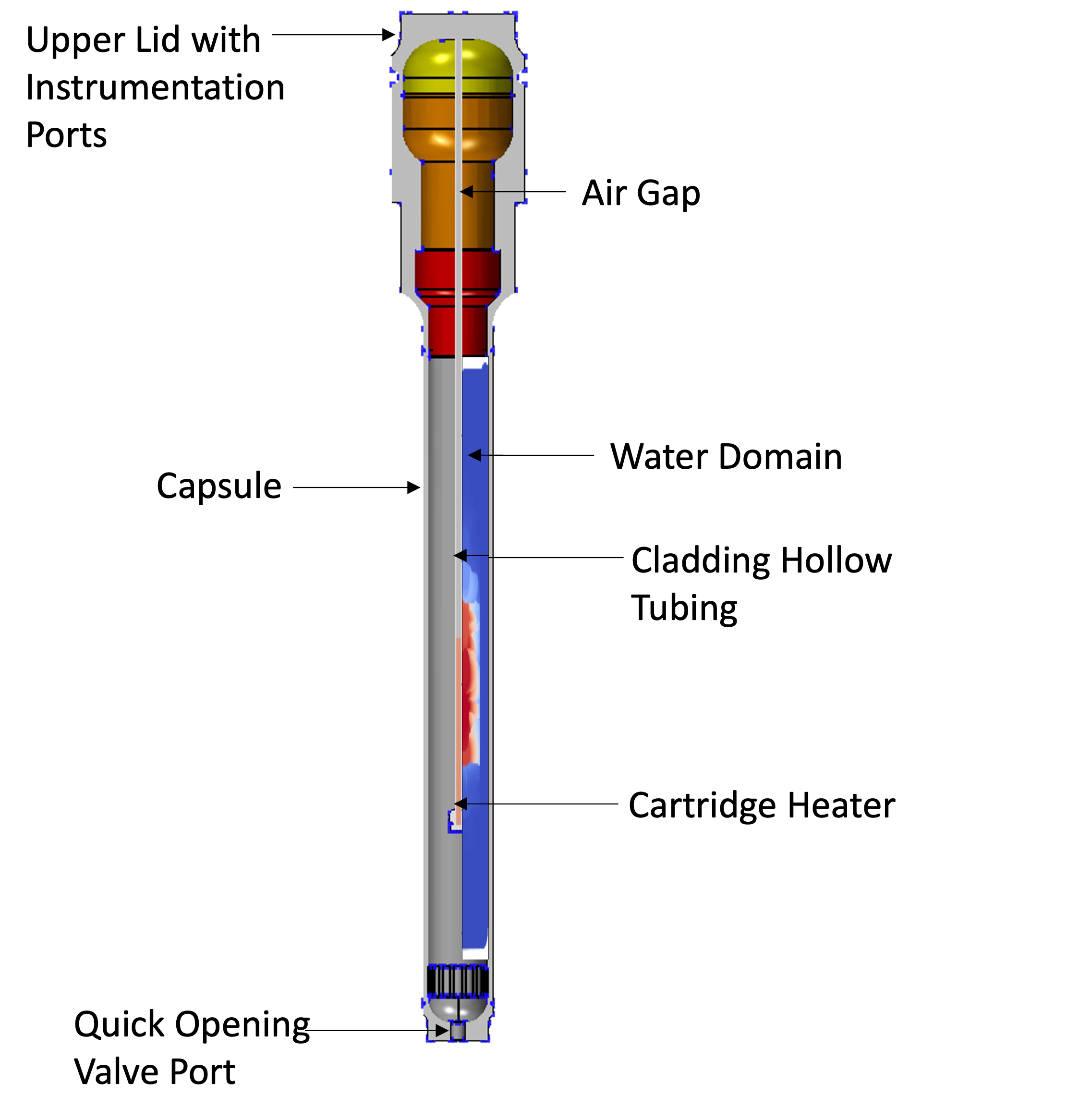}
		\caption{Schematic}
	\end{subfigure}
	\begin{subfigure}{.03\textwidth}
     \begin{overpic}[scale=.21,percent]{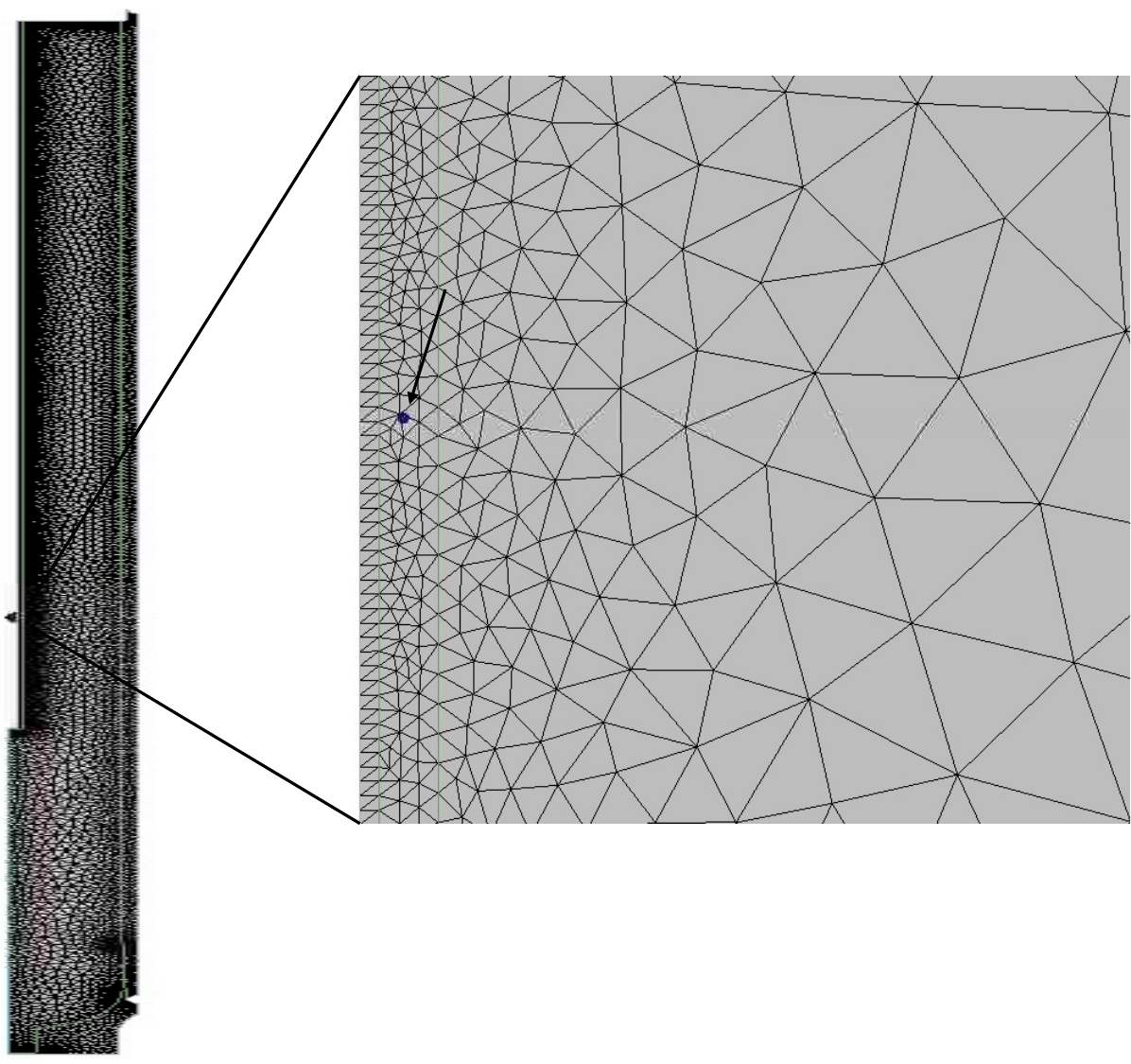}
          \put(40,73){Heat Source}
        \put(40,43){Heater Region}
    \end{overpic}
         \caption{\label{fig:TWIST_mesh}\footnotesize{Mesh}}
	\end{subfigure}\hspace{2em}

\caption{\label{fig:TWIST_geometry}\textbf{Geometry and mesh of OPTI-TWIST.}  
The axial symmetry of the OPTI-TWIST is exploited by simulating only half the domain as the cartridge heater is placed at the center of the capsule. The geometry and mesh reveal richer dynamics near the heater region. }
\end{figure}

Here, our sensor placement optimization is demonstrated on the OPTI-TWIST prototype, which is electronically heated to mimic the neutronics effect of TWIST prior to insertion into a reactor at 
INL. 
Temperature is the field of interest, and point thermocouples will be used as the sensors. 
The OPTI-TWIST prototype was designed to simulate thermal-hydraulics behavior of TWIST during irradiation in the reactor, as well as to measure the effect of loss of coolant on the fuel rodlet temperature. 
In OPTI-TWIST, the fuel-rod specimen is replaced by an instrumented electric cartridge heater, and loss of coolant is controlled by a quick-opening valve at the bottom of the capsule.  
To provide the temperature fields necessary to train the sparse sensing algorithm, a simplified 2D CFD model of the OPTI-TWIST geometry is developed using StarCCM+ (\autoref{fig:TWIST_geometry})~\cite{ccm0}. 

The CFD model accounts for steady-state turbulent natural circulation conditions, including two controlled parameters: heater power ($\dot{q}$) and outer surface temperature ($T_{sur}$). 
\begin{figure*}[hbt!]
\begin{subfigure}{.49\textwidth}
    \includegraphics[width=\textwidth]{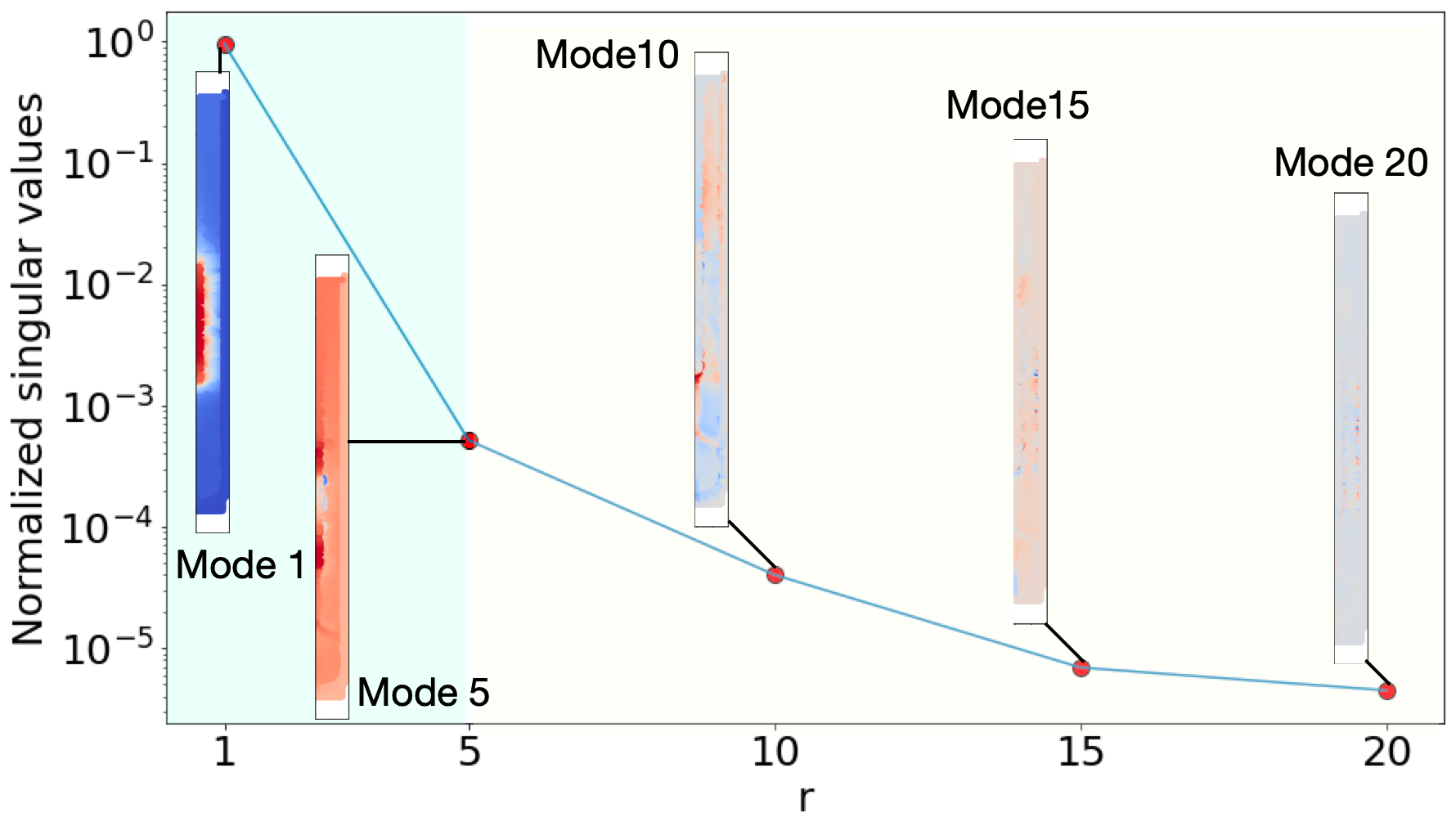}
    \caption{\label{fig:TWISTPOD}Leading POD modes and energy contribution.}
\end{subfigure}
\begin{subfigure}{.13\textwidth}
    \begin{overpic}[scale=.29,percent]{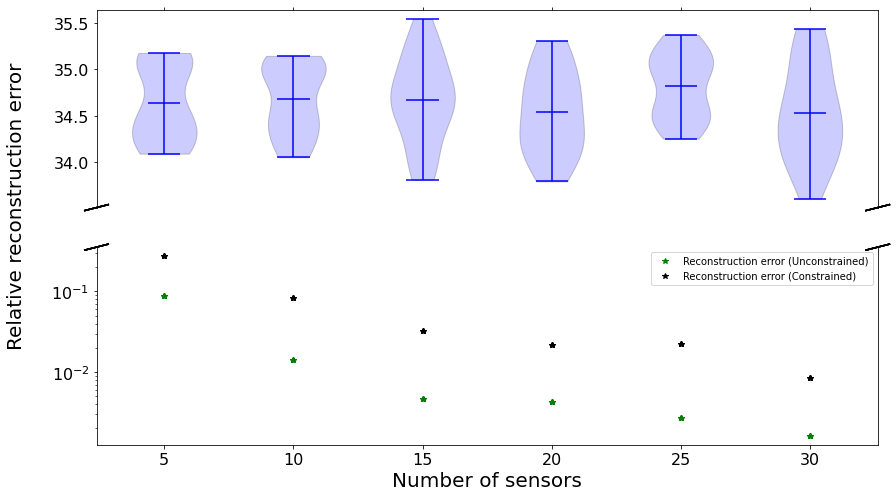}
\put(65,2){$(r)$}
  \end{overpic}
    \caption{\label{fig:TWIST_errors} Error Study}
\end{subfigure}
	\caption{\label{fig:TWIST_POD_errors} \textbf{An informed trade-off between reconstruction accuracy and number of sensors is possible for QR optimized sensors.} The leading POD modes capture 99\% of energy content and just 5 sensors are enough to obtain a accurate reconstruction with $\epsilon \sim \mathcal{O}(10^{-1})$ (a). QR selected sensor accuracy increases with an increase in the number of sensors as compared to random placements that produce orders of magnitude larger relative reconstruction errors (b). }
\end{figure*}
These two controlled parameters (i.e., heater power and surface temperature) were varied, while keeping the initial temperatures ($T_0=300$K) and the system pressure ($P_{sys}=2250$psi) constant throughout. The data are comprised of 49 steady-state temperature fields resulting from seven heater powers and surface temperatures uniformly sampled at 350--650W and 240--420K, respectively. The convergence criterion was the maximum liquid temperature, which showed negligible fluctuations after 2000 time steps. 
\begin{table}[b!]
\small
\begin{center}
\begin{tabular}{|c| c|c|} 

 \hline
 \footnotesize{Sensor Placement} & \footnotesize{Optimization Objective} & \footnotesize{Reconstruction Error}($\epsilon$)\\ [0.5ex] 
 \hline\hline
 \footnotesize{Unconstrained}  & \footnotesize{5.432829081027846e-10} & \footnotesize{0.168}\\ 
 \hline
\footnotesize{Constrained} & \footnotesize{4.534195929074569e-11} & \footnotesize{0.174} \\
 \hline
 \footnotesize{Random}  &  \footnotesize{1.026196077627373e-12} & \footnotesize{25.24} \\

 \hline

\end{tabular}
\caption{\label{rec_error}Summary of the relative reconstruction error ($\epsilon$) and optimization criteria ($\log |\det\Sel\boldsymbol{\Psi_r}|$) for sensor placement given in \autoref{fig:TWISTprofile}.}
\end{center}
\end{table}
First, our optimization is run on the steady-state temperature fields, resulting in the unconstrained optimal placement shown in \autoref{fig:TWIST1}. Unconstrained optimization selects three sensors near the heater (\autoref{fig:TWIST1}); however, space restrictions make these heater-adjacent locations experimentally infeasible. Enforcing sensor constraints to lie outside the heater region results in a reconstruction error of $\epsilon= 0.174$---an increase of only .006. \autoref{fig:TWISTrand} contrasts these optimized sensor reconstructions with ensembles of randomly placed sensors. Observe that $\epsilon_{unconstrained} < \epsilon_{constrained} \ll \epsilon_{random}$, i.e. placing sensors in random locations leads to significantly larger reconstruction errors (\autoref{rec_error}).

Gaussian i.i.d. noise $\boldsymbol{\eta} \sim \mathcal{N} (0,0.01)$ is added to the measurements to analyze the uncertainty heatmaps in reconstruction of the true temperature profile. Removing sensors from heater-adjacent locations leads to an increase of approximately 0.5K in the uncertainty in reconstruction of the flow field close to the heater (\autoref{fig:TWIST_unc_const}) compared to uncertainty resulting from unconstrained sensor placement (\autoref{fig:TWIST_unc_unconst}). The uncertainty in reconstruction is higher by 40-50K throughout the domain when sensors are placed randomly (\autoref{fig:TWIST_unc_rand}). Stratified contours of reconstruction errors occur where random sensors fail to accurately capture the transition between hot and cold (\autoref{fig:TWISTrand}). Therefore, randomly placed sensors fail to capture heater-correlated fluctuations and result in higher reconstruction errors and uncertainty.

\,The leading two POD modes, which drive approximately 99\% of the energy content, capture the heat advection from the heat source. Thus, only five sensors---corresponding to the first five POD modes---are required to reconstruct the flow with high accuracy. The cumulative energy content, along with a visualization of the first three POD modes, is given in \autoref{fig:TWISTPOD}. These POD modes capture the interfaces between lower and higher temperatures as the advection flow progresses for different operating ranges.
Sensors placed at random locations fail to capture the underlying physics of the system. Increasing the number of random sensors selected does not significantly improve the reconstruction of the flow field. An ensemble of sensors placed at random locations produces large relative reconstruction errors that average approximately 35 as the number of sensors is increased from 5 to 30, as seen in \autoref{fig:TWIST_errors}.
Random placement of sensors increases the reconstruction error by ten orders of magnitude. The random placement strategy alludes to data-agnostic sensor placement at best-guess locations for reconstructing temperature fields. 

Unconstrained optimization favors locating sensors close to the heat source, due to the richer dynamics that exist there. Imposing sensor constraints within QR results in a near-optimal placement outside this region, as well as negligible loss of reconstruction performance. 
Moreover, the error decreases with more optimized sensors (unconstrained or constrained); however, random placements still suffer from high error even with additional sensors.
Therefore, placing or adding sensors without optimizing them in regard to the underlying flow or deployment constraints can introduce large errors in the corresponding digital twins especially in the presence of noisy sensor measurements---and may even cause sensor damage.
Incorporating such considerations prior to setting up a physical experiment enables precise uncertainty quantification and helps validate the predictions of a digital twin.

\vspace{-1em}
\subsection{Transient simulation of the OPTI-TWIST prototype} 
\label{ssec:OPTI-TWIST_Transients}

\begin{figure*}[t]
\begin{subfigure}{.49\textwidth}
    \includegraphics[width=\textwidth]{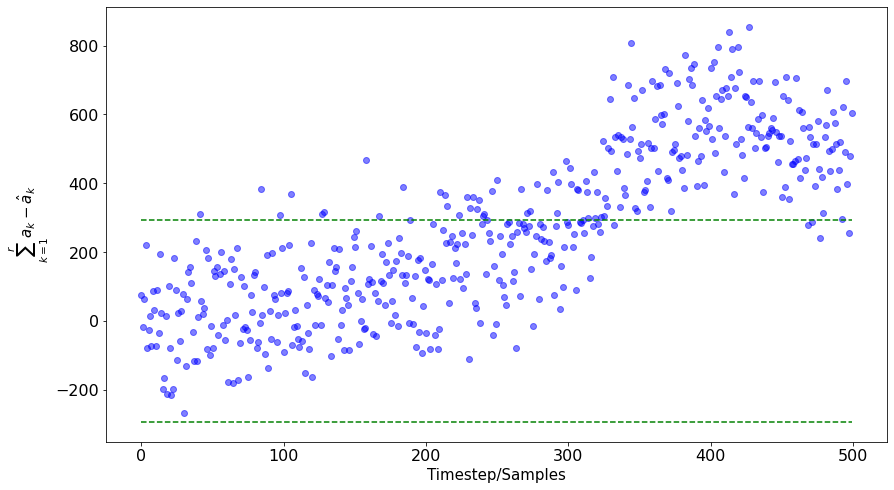}
    \caption{\label{fig:TWIST_transient1}\footnotesize{Rank 10 model}}
\end{subfigure}
\begin{subfigure}{.49\textwidth}
    \includegraphics[width=\textwidth]{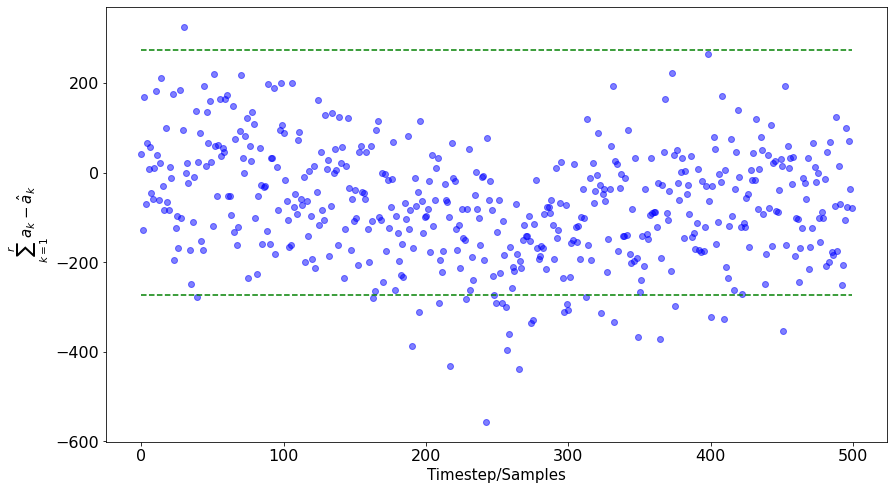}
    \caption{\label{fig:TWIST_transient2}\footnotesize{Rank 20 model}}
\end{subfigure}
\begin{subfigure}{.49\textwidth}
    \includegraphics[width=\textwidth]{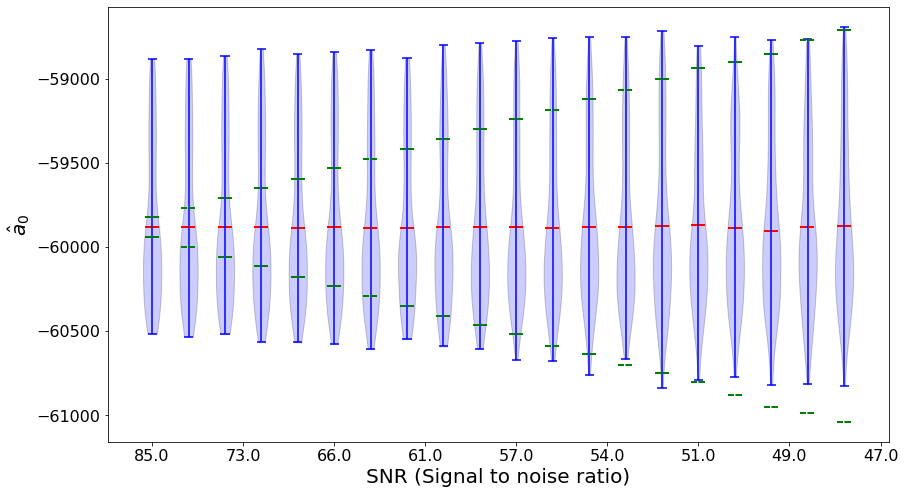}
    \caption{\hspace{-0.81em}\label{fig:TWIST_transient3} \footnotesize{Increasing uncertainty in estimation of $\hat{\ba}_k$ with a decreasing SNR.}}
\end{subfigure}
\begin{subfigure}{.49\textwidth}
    \includegraphics[width=\textwidth]{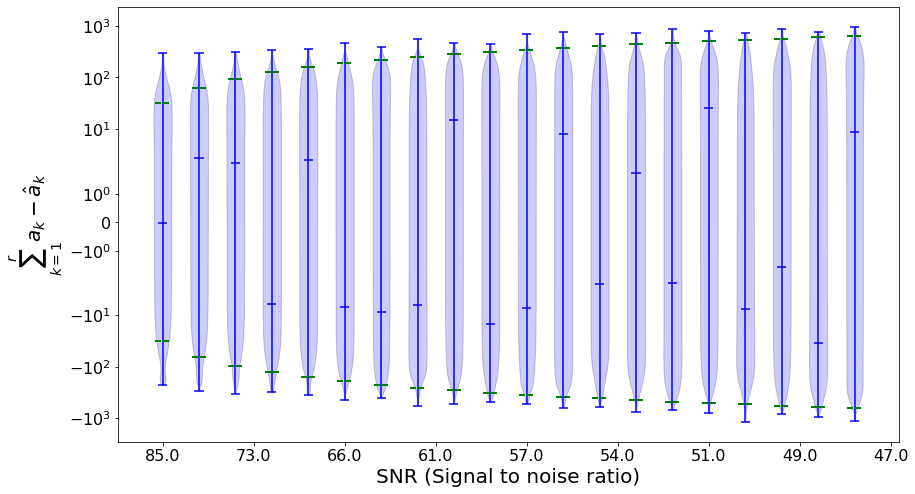}
    \caption{\label{fig:TWIST_transient4}\footnotesize{Increase in estimation errors with a decreasing SNR.}}
\end{subfigure}
	\caption{\label{fig:TWIST_POD_errors} \textbf{Estimation error analysis for OPTI-TWIST heater transients.} (a) Uncertainty estimation reveals that a rank 10 model is not sufficiently descriptive of dynamics after t =300 under sensor noise. The rank 20 approximation (b) is valid over a longer time horizon of 500s of test data. As the SNR decreases, (c) POD coefficient variance increases, which propagates to an increase in uncertainty of estimation errors (d). }
\end{figure*}

\begin{figure}[t!]
    \begin{subfigure}{.13\textwidth}
        \begin{overpic}[scale=.47,percent]{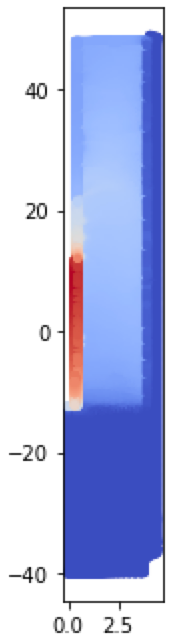}
          \put(-2.5,34){\rotatebox{90}{Height [cm]}}
    \end{overpic}
		
		\caption{\hspace{-0.7em}\label{fig:TWISTgroundTruthTransient} \scriptsize{Ground Truth}}
	\end{subfigure}
 \hspace{0.5em}
	\begin{subfigure}{.25\linewidth}
            \includegraphics[width=0.65\textwidth]{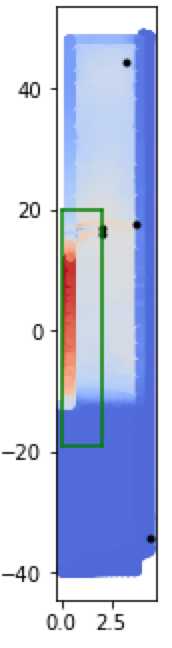}
		\caption{\label{fig:TWISTconstrainedtransient}\scriptsize{Constraint}}
	 \end{subfigure}
  \hspace{1.7em}
	\begin{subfigure}{.155\textwidth}
		\includegraphics[width=\textwidth]{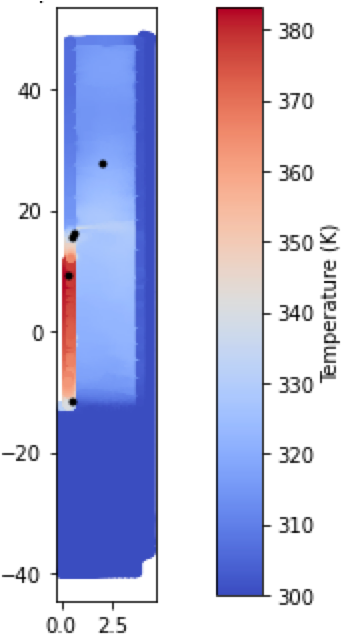}
		\caption{\label{fig:TWISTQRTransient}\scriptsize{QR}}
	\end{subfigure}

\caption{\textbf{Sensor placements superimposed on reconstructions for heater power transients.} Heater-adjacent temperature fluctuations result in (c) sensors optimized close to the heater and a corresponding low reconstruction error $\epsilon$ = 1.026. When constraints are imposed, sensors are placed outside the green constraint region (b), resulting in higher reconstruction error $\epsilon$ = 2.042.}
\end{figure}

Analyzing the effect of power transients on reactor core coolant temperature, pressure, and velocity is essential for real-time safety monitoring and control of a nuclear reactor. Here, we optimize sensor placements to capture the dynamics of the heat flow when the heater power is varied as a transient. During transients it is essential to capture the instance when sensor readings start diverging from predicted metrics in the presence of noise. This can signal the need for model recalibration and can prevent accidents caused by power surges at a nuclear facility. The data is obtained from the 2D CFD model described in \autoref{ssec:OPTI-TWIST} which runs for 600s, and is comprised of 1000 temperature fields sampled at every 0.6s. The heater transient power profile $P(t)$ can be described by

\begin{equation}
\begin{aligned}
  \label{eqn:heater_power_transient}
    P(t) &=  P_o,\, \text{if} \, t \leq t_1, \\
    P(t) &= P_o + \frac{P_2 - P_o}{t_2 - t_1} * (t-t_1), \: \text{if} \: t_1 \leq t \leq t_2, \\ 
    P(t) &= P_2, \, \text{if} \, t_2 \leq t \leq t_3,
\end{aligned}
\end{equation}
where $P_o = 10 W, P_2 = 250 W, t_1 = 200 s, t_2 = 2t_1 \text{and}, t_3 = 3t_1$. Similar to the steady state temperature profile, richer dynamics are located in heater adjacent regions (\autoref{fig:TWISTgroundTruthTransient}). The unconstrained sensor layout shows a number of sensors near the heater (\autoref{fig:TWISTQRTransient}), however due to spatial constraints all sensors must be located away from the heater (\autoref{fig:TWISTconstrainedtransient}). The increase in reconstruction error due to imposing constraints is as low as 1\%. The algorithm is trained on the first 500 timesteps and reconstructs the temperature profiles at the last 500 timesteps from optimized sensors readings with additive i.i.d . Gaussian noise $\boldsymbol{\eta} \sim \mathcal{N} (0,0.9)$.

We study the reconstruction uncertainty using the predicted distribution of estimation coefficients in \autoref{sec:UQ}. As shown in \autoref{fig:TWIST_POD_errors}, the rank 10 gappy POD model is insufficient to characterize the transient behavior in test data and reconstruction errors increase beyond the established bounds at $t=200s$, whereas the rank 20 model captures the dynamics more effectively over the entire time horizon. When information regarding the true state is available, model recalibration can be signaled by bounds established for the error distribution in the presence of increasing noise (\autoref{fig:TWIST_transient4}). When true coefficients are unavailable, the bounds established for the distribution of estimated POD coefficients $\hat{\ba}_i$ can be used to flag erroneous readings. This uncertainty estimation of the flow field during transients ensures safe operating conditions during experimentation and testing of nuclear reactors.

\section{Discussion and Outlook}

The reconstruction of reactor core flow fields using a limited budget of sensors has emerged as a critical enabler for digital twinning of nuclear assets. However, achieving optimal sensing in high-dimensional real-world systems extends beyond the nuclear industry and encompasses diverse fields such as biology, physics, aviation, and automotive industries. Engineering systems are subject to different constraints and limitations on sensors, making the selection and optimal placement of sensors while considering these constraints a crucial aspect of algorithm development. In this study, we have developed strategies for placing sensors to satisfy constraints related to sensor proximity, regional limitations on sensor quantity, and design or user prescribed locations. Through these strategies, we have demonstrated the effectiveness of adaptive sensor placement in satisfying constraints while maintaining optimality and accuracy of the reconstructed responses of interest. Moreover, we showcase the scalability and broad applicability of the algorithms on a variety of applications and constraints.

 Nevertheless, more complex constraints may arise in nuclear, fluid, or aerospace applications in which the capability to achieve flow reconstruction based on sparse measurements is indispensable. For instance, in each reactor region the maximum number of allowable sensors is usually design-specific and cannot be exceeded. Moreover, an emerging practice during nuclear fuels tests is the use of distributed sensors (e.g., fiber-optic sensors or multipoint thermocouples). In fiber bragg grating (FBG), the refraction index changes along the sensor length and provides distributed measurements. Designing fiber optic sensors---which act as line sensors with different measurements at each point---is a novel challenge, as line sensors can be topologically shaped along various structures in engineering systems. Optimizing such topologies is another future direction for adaptive sensor placement. Fiber optic bundles are used for recreating high-quality images in both nuclear engineering and neuroscience. Optimizing the locations for these bundles to capture the best quality images is another interesting research direction.

Furthermore, it might be very costly to place sensors in certain areas of the reactor, due to the need for specially designed sensors capable of withstanding harsh working conditions. Other areas may entail spatial constraints. Thus, multi-objective optimization based on optimizing the sensor cost, spatial locations, as well as predicted dynamics will be essential. Time-dependent dynamics and the study of transients is invaluable in the nuclear field. Sensor placement based on time-dependent data from OPTI-TWIST and the use of dynamic mode decomposition or a nonlinear embedding such as autoencoders can help generalize to new physics scenarios, and will require new uncertainty estimates that can handle bias inherent to these types of models.

The ultimate goal is to extend the algorithms to inform users of optimal locations and timesteps to collect spatiotemporal sparse measurements to reconstruct core flow profiles during power transients. This should naturally evolve to the capability of performing optimal sensor placement for multi-class classification, where the algorithm must select the sensor network capable of predicting which accident scenario is more likely to occur faster than real time. Examples of such accident scenarios include Loss of Coolant Accident (LOCA), Reactor Initiated Accidents (RIA), and loss of power. These scenarios can be easily realized in a non-destructive fashion within the TWIST prototype by opening a valve, or suddenly shutting the heater power off. Data-driven sensing frameworks have the potential to identify sensor maps capable of detecting off-normal conditions, anomalies, and injected signals, enhancing resilience and security of the physical twin against cyber-attacks.
\section*{Acknowledgments}
The material presented herein was based on work supported through the Idaho National Laboratory (INL) Laboratory Directed Research \& Development (LDRD) Program under DOE Idaho Operations Office Contract DE-AC07-05ID14517 for LDRD-22A1059-091FP. The material presented here is implemented in PySensors \cite{deSilva2021}, an open source, Scikit-learn style Python package for the sparse placement of sensors and in RAVEN (Risk Analysis and Virtual ENviroment), INL's in-house open source software \cite{rabiti2021raven}. K. Manohar and S. Brunton acknowledge support from the NSF AI Institute in Dynamic Systems under grant 2112085. The authors would like to thank Andrei A. Klishin for helpful discussions.


\printbibliography

\end{document}